\documentclass[12pt]{article}
\usepackage{color}
 \usepackage{amscd,amsfonts,amssymb,amscd,amsmath,latexsym}
\usepackage[hypertex]{hyperref}
\usepackage[all]{xy}
\usepackage{mathrsfs}
\usepackage{graphicx}

\title{Secondary Invariants for String Bordism and tmf}
\author{Ulrich Bunke\thanks{NWF I - Mathematik,
Universit{\"a}t Regensburg,
93040 Regensburg,
GERMANY, ulrich.bunke@mathematik.uni-regensburg.de}, Niko Naumann \thanks{NWF I - Mathematik,
Universit{\"a}t Regensburg,
93040 Regensburg,
GERMANY, niko.naumann@mathematik.uni-regensburg.de} }

\newtheorem{theorem}{Theorem}[section] 
\newtheorem{prop}[theorem]{Proposition}
\newtheorem{lem}[theorem]{Lemma}
\newtheorem{ddd}[theorem]{Definition}
\newtheorem{kor}[theorem]{Corollary}

\newcommand{\modu}{{\tt mod\:}}
\newcommand{\Filt}{{\tt Filt}}

\newcommand{\F}{{\mathbb{F}}}

\renewcommand{\P}{{\mathbb{P}}}

\newcommand{\Z}{\mathbb{Z}}

\newcommand{\diag}{{\tt diag}}
\newcommand{\proof}{{\it Proof.$\:\:\:\:$}}
\newcommand{\ori}{{\tt or}}

\newcommand{\R}{\mathbb{R}}

\newcommand{\Q}{\mathbb{Q}}

\renewcommand{\sinh}{\mathrm{sinh}}

\newcommand{\C}{\mathbb{C}}

\newcommand{\Tr}{{\tt Tr}}

\newcommand{\cZ}{\mathcal{Z}}

\newcommand{\cM}{\mathcal{M}}

\newcommand{\cC}{\mathcal{C}}

\newcommand{\cO}{\mathcal{O}}

\newcommand{\Ext}{{\tt Ext}}

\newcommand{\im}{{\tt im}}
\newcommand{\sign}{{\tt sign}}

\newcommand{\id}{{\tt id}}
\newcommand{\ord}{{\tt ord}}
\newcommand{\nat}{\mathbb{N}}

\newcommand{\spec}{{\tt spec}}

\def\imath{{i}}

\def\hB{\hspace*{\fill}$\Box$ \newline\noindent}

\newcommand{\ind}{{\tt index}}

\def\hB{\hspace*{\fill}$\Box$ \\[0.5cm]\noindent}

\newcommand{\ch}{{\mathbf{ch}}}

\newcommand{\bV}{\mathbf{V}}

\newcommand{\hA}{\hat{\mathbf{A}}}

\begin{document}\maketitle\begin{abstract}
Using spectral invariants of Dirac operators we construct a secondary
version of the Witten genus, a bordism invariant of string manifolds in dimensions $4m-1$.
We prove a secondary index theorem which relates
this global-analytic construction with its homotopy-theoretic analog. The latter will be calculated
through its factorization over  topological modular forms.
 \end{abstract}

\tableofcontents
 
\section{Introduction}

\subsection{Review of the contents}

A well-established principle in the border area between algebraic topology and global analysis predicts that local geometric representations of primary topological invariants lead to 
interesting secondary constructions \cite[Introduction]{bunke-2008}.
In the present paper, based on this principle,  we  construct an interesting secondary invariant for
string bordism classes  and elements of the homotopy of the spectrum of  topological modular forms.
We will employ the notion of geometric string 
structures which was  recently developed by Waldorf \cite{waldorf}. 

Here is a short summary of our constructions. For all $m\ge 1$ we define an abelian group
$T_{2m}$ (Definition \ref{uzifefwfwefewfwef}) and construct a homomorphism (Definition \ref{uzduiwqdqwd})
\begin{equation}\label{bbb}
b^{an} :MString_{4m-1}\to T_{2m}
\end{equation}
using spectral invariants of Dirac operators  and differential geometry. 
It is a secondary version of the Witten genus.
Because of the appearance of spectral invariants, a direct evaluation of $b^{an}$ is complicated. But it can be analyzed on the subgroup
\begin{equation}\label{bbb333}A_{4m-1}:=\ker(MString_{4m-1}\stackrel{j}{\to} MSpin_{4m-1})\subseteq MString_{4m-1}\ .\end{equation}
The restriction of $b^{an}$ to $A_{4m-1}$ will be compared with homomorphisms
\begin{equation}\label{bbb1}
b^{geom},b^{top} :A_{4m-1}\to T_{2m}
\end{equation}
defined using  differential geometry (\ref{w8e9fwefwfwefwf}), and  homotopy theory (Definition \ref{uzidqwdqwdwqdwd}), respectively.
 The constructions of the homomorphisms  $b^{an}$ and $b^{geom}$   are very similar to the analytic and geometric  description of the Adams $e$-invariant given in \cite{MR0397797}. Our first main result asserts that 
\begin{equation}\label{uiwefwefwefwef6565}
b^{an}_{|A_{4m-1}}=b^{geom}=b^{top}\ .
\end{equation}
The first equality $b^{an}_{|A_{4m-1}}=b^{geom}$ is essentially a consequence of the Atiyah-Patodi-Singer index theorem \cite{MR0397797} and will be verified during the construction of $b^{geom}$ in Subsection  \ref{ioqdqwd}. The second equality
$b^{geom}=b^{top}$, shown in Theorem \ref{e2e7uiwqdqwdwqd}, is deeper, and the resulting equality $b^{an}_{|A_{4m-1}}=b^{top}$ can be considered
as a secondary index theorem. In the proof we use  ideas  we  learned from  \cite{laures-cob}.
The analogy with the secondary index theorem for the $e$-invariant of Adams \cite{MR0159336} as formulated  in \cite[Introduction]{bunke-2008} will be explained in Subsection \ref{uidqwdqwdqwdqwdee}.

 We will give a complete and explicit homotopy theoretic calculation of $b^{top}$. It is based on the factorization (Proposition \ref{udiqwdqwdqwdq})
$$\xymatrix{MString_{4m-1}\ar[r]^\sigma&tmf_{4m-1}\ar[d]^{b^{tmf}}&\\A_{4m-1}\ar@{_{(}->}[u]\ar[r]^{b^{top}}&T_{2m}}
$$ and the explicit calculation of $b^{tmf}$ (Propositions  \ref{uidodqwdqwd} and \ref{zwdqdwqd}) whose construction is similar to the one of  $b^{top}$.
Here, $tmf$ denotes the connective spectrum of topological modular forms of Goerss-Hopkins-Miller and Lurie \cite{goerss1}, and $\sigma$ is induced by the $tmf$-valued Witten genus which has been constructed as an $E_\infty$-ring map $\sigma :MString\to tmf$  in \cite{ahr}.

Since the homotopy groups of string bordism
and $tmf$ are interesting and complicated, it is useful to relate these objects non-trivially to geometry or global analysis.

\bigskip
{\em Acknowledgements: The first author is very grateful to K. Waldorf for introducing  him to the notion  of a geometric string structure
and further discussion at earlier stages of this work. He further thanks
M. Hovey for very helpful remarks about the structure of the string bordism groups. 

The charts of the spectral sequences have been compiled using the TEX-package of T. Bauer.
The homotopy groups of $MSpin$ have been calculated with the help of  MAPLE.}

\subsection{Thom spectra}

In the remainder of this introduction we give a more detailed overview.
We first recall the definition of the Thom spectra appearing in the present paper.
Let \begin{equation}\label{uifwefwefwefewf}
\dots \to  BO\langle 8 \rangle\to BO\langle 4\rangle\to BO\langle 2\rangle\to BO
\end{equation}
be the first stages of the Postnikov tower of the $H$-space $BO$. The classical names of these spaces are
$$BSO=BO\langle 2\rangle\ ,\quad BSpin=BO\langle 4\rangle\ ,\quad BString = BO\langle 8 \rangle\ .$$

There is a functor from the category of $H$-spaces over $BO$ to ring spectra which associates to a map $\xi:X\to BO$ the  Thom spectrum $X^\xi$.
The Thom spectrum for the map $BO\langle i \rangle\to BO$ will be denoted with $MO\langle i\rangle$. In particular we will write
$$MSO=MO\langle 2\rangle\ ,\quad MSpin=MO\langle 4\rangle\ ,\quad MString = MO\langle 8 \rangle\ .$$ These spectra are related by morphisms of ring spectra
$$MString \stackrel{j}{\to} MSpin \to MSO\to MO\ .$$

\subsection{The Witten genus, $T_{2m}$, and $b^{geom}$} 

For the purpose of this introduction, we explain the construction of $b^{geom}$ which is the simplest of the three maps in (\ref{uiwefwefwefwef6565}).  In particular we want to motivate the Definition \ref{uzifefwfwefewfwef} of the group $T_{2m}$.
The map $b^{geom}$ is a secondary version of    the  Witten genus.
We let $KO$ denote the real $K$-theory spectrum. For a (ring) spectrum $X$ we let
$X[[q]]$ denote the (ring) spectrum which represents the homology theory
$X[[q]]_*(B):=X_*B[[q]]$ for all finite $CW$-complexes $B$.
The Witten genus is a morphism of ring spectra 
$$R:MSpin\to KO[[q]]\ .$$ 
We refer to  Subsection \ref{udqwidqwdqwdqwd444} for a description of this map as a transformation of homology theories. 

For a spectrum $X$ we let $X_* :=\pi_*(X)$ denote the homotopy groups of $X$.  
Via the Thom-Pontrjagin construction, elements in the homotopy of  Thom spectra like $MSpin$ or $MString$ can be realized as bordism classes of manifolds with corresponding structures on the tangent bundles (see Subsection \ref{udidqwdqwdqwdqwdqwdqwd}).
Thus homotopy classes of $MString$ correspond to bordism classes of closed string manifolds $[M,\alpha^{top}]$, where a string manifold $(M,\alpha^{top})$ consists of a spin manifold $M$ together with the choice of a topological string structure $\alpha^{top}$, see Subsection \ref{udidqwdqwdqwdqwdqwdqwd} for details. The morphism $j:MString_*\to MSpin_*$ just forgets the topological string structure, i.e. it is given geometrically by $j([M,\alpha^{top}])=[M]$.

We now give a description of the Witten genus following the exposition \cite{MR1189136}.
Let $[M]\in MSpin_{4m}$ be a homotopy class represented by a closed $4m$-dimensional spin manifold $M$.  We choose a connection $\nabla^{TM}$ on the tangent bundle $TM$ of $M$. In general,
if $\nabla^V$ is a connection on a real vector bundle $V\to B$ over some manifold $B$, then by $p_i(\nabla^{V})\in \Omega^{4i}(B)$ we denote the Chern-Weil representative of the $i$'th
Pontrjagin class $p_i(V)\in H^{4i}(B;\Z)$. We set $\kappa_{m}:=1$ for even $m$, and $\kappa_{m}:=\frac{1}{2}$ in the case of odd $m$.
In (\ref{udiqwdqwdqwd3}) we define a power series
$$\Phi\in \Q[[q]][[p_1,p_2,\dots]]\ .$$  With this notation
we have the following local expression for the Witten genus
\begin{equation}\label{uiqdqwdqwdwqdwqd}
R([M])=\kappa_{m}
 \int_M  \Phi(p_1(\nabla^{TM}),p_2(\nabla^{TM}),\dots)\ ,
\end{equation}
which a priori is an element of $\Q [[q]]$.
By the Atiyah-Singer index theorem, the coefficients of this formal power series
can be identified with  indices of twisted Dirac operators. More precisely, we have 
$$R([M])=\kappa_m\sum_{n\ge 0} q^n\ind(D_M\otimes R_n(TM))\ ,$$
where   $R_n(TM)$  is a certain virtual bundle derived from the tangent bundle (see (\ref{zfquefhefwefwef})) for all $n\ge 0$, and $D_M\otimes  R_n(TM)$ denotes the spin Dirac operator of $M$ twisted by 
$R_n(TM)$.  Since the indices are integers and even in the case of odd $m$ (because of  an additional real symmetry) we see that
$$R([M])\in \Z[[q]]\ .$$

 For a subring $R\subseteq \C$  we denote by $\cM^R_{2m}$ the space of modular forms for $SL(2,\Z)$ of weight $2m$   with a  $q$-expansion in the subring $R[[q]]\subseteq \C[[q]]$.
Using the $q$-expansion   we will identify $\cM^R_{2m}$ with a sub-$R$-module  of $R[[q]]$.
For more information on modular forms we refer to Subsection \ref{fzuefwefwefewf}.

It is a crucial observation for our constructions that the Witten genus   has the following factorization
$$\xymatrix{&&&\cM^\Z_{2m}\ar@{_{(}->}[d]\\MString_{4m}\ar@{.>}[urrr]\ar[r]_j&MSpin_{4m}\ar[r]_R&KO[[q]]_{4m}\ar[r]_\cong&\Z[[q]]}\ .$$

The construction of the secondary  invariant $b^{geom}$  starts with the local formula (an integral over characteristic forms)  for the Witten genus given by the right-hand side of (\ref{uiqdqwdqwdwqdwqd}). It will be applied to a spin zero-bordism $Z$ of a string manifold $(M,\alpha^{top})$ and gives, combined with a contribution from a geometric  string structure $\alpha$, a formal power series in $\R[[q]]$.
In order to ensure independence of the choices of the geometry and the zero bordism we have to  calculate modulo everything that comes from closed spin manifolds
and string zero bordisms.  
This leads us to the definition of the following quotient group.
\begin{ddd}\label{uzifefwfwefewfwef}
The group $T_{2m}$  
is defined by 
$$T_{2m}:=\frac{\R[[q]]}{\Z[[q]]+\cM_{2m}^\R}\ .$$
\end{ddd}

We now turn to the construction of the secondary invariant $b^{geom}(M,\alpha^{top})$.
It is based on the existence of a spin zero bordism of the string manifold $(M,\alpha^{top})$. We therefore have to  assume that
the class $[M,\alpha^{top}]\in MString_{4m-1}$ is in the subgroup
$A_{4m-1}\subseteq MString_{4m-1}$ defined in (\ref{bbb333}).

Let  $(M,\alpha^{top})$ be a $4m-1$-dimensional string manifold such that $[M,\alpha^{top}]\in A_{4m-1}$.
Then we can choose     a spin zero bordism $Z$ of $M$. Furthermore, we choose a  connection $\nabla^{TZ}$ extending a connection $\nabla^{TM}$ on the tangent bundle $TM$ of $M$. 

A topological string structure $\alpha^{top}$ is by definition a "trivialization" of the spin characteristic class $\frac{p_1}{2}(TM)\in H^4(M;\Z)$. A geometric refinement $\alpha$ of  $\alpha^{top}$ trivializes
this class on the form level, i.e. $\alpha$ gives rise to a form
$H_\alpha\in \Omega^{3}(M)$ such that $dH_\alpha=\frac{1}{2} p_1(\nabla^{TM})$
(see \cite{waldorf}).
This can be used to define a refinement
$\tilde p_1(\nabla^{TZ},\alpha)\in \Omega^4(Z)$ of $p_1(\nabla^{TZ})$
which vanishes on the boundary $M$ of $Z$ (see Subsection \ref{udqwidqwdqwdqwd}).
The value  of $b^{geom}$ on the class $[M,\alpha^{top}]\in MString_{4m-1}$ is then given  (compare Lemma \ref{zdqwud}) as the class
$$b^{geom}([M,\alpha^{top}]):=\left[\int_Z \Phi(\tilde p_1(\nabla^{TZ},\alpha),p_2(\nabla^{TZ}),\dots)\right]\in T_{2m}\ .$$

In contrast to $b^{geom}$, the analytic variant $b^{an}([M,\alpha^{top}])\in T_{2m}$
 in  Definition  \ref{uiqdwqdqwdwqd}     provides an intrinsic formula which does not depend on the choice of a spin zero bordism and is therefore defined on all of $MString_{4m-1}$. It  involves spectral invariants of  the twisted Dirac operators $D_M\otimes R_n(TM)$ for all $n\ge 0$.
 In Lemma \ref{udidqwdqwd} we show the equality
$$b^{geom}=b^{an}_{|A_{4m-1}}\ .$$

 The construction of all variants of $b$ is based on the interplay between the facts that $MString_{4m-1}$ is torsion and $KO[[q]]_{4m}\otimes \Q$ is non-trivial. This explains the restriction   to dimensions of the form $4m-1$.
In order to detect elements of $MString_*$ in dimensions $0,1,2,4 \:\:\mbox{mod $8$}$ one can use the Witten genus $R\circ j:MString\to KO[[q]]$ directly.

\subsection{Calculations}

In Section \ref{gfdhfdf} we consider the case of three-manifolds, i.e. the case  $m=1$.
In this case it suffices to consider the constant term of the formal power series
representing $b^{geom}([M,\alpha^{top}])$. This simplifies matters considerably and justifies
a separate discussion. We will see that
$$b^{geom}:MString_{3}\cong \Z/24 \Z\hookrightarrow T_2=\R[[q]]/\Z[[q]]$$ is injective.
As a side result we get the following analog for spin manifolds of Atiyah's canonical $2$-framings of oriented three-manifolds \cite{MR1046621}: {\em A $3$-dimensional connected closed spin manifold has  a canonical
topological string structure (Definition \ref{uiqdwqdqwdw})}.

  In higher dimensions the homotopy groups $MString_*$ are not fully understood. The following facts are known \cite{MR1455523}.
\begin{enumerate}
\item $MString_{k}$ is a finite group for $k\equiv 1,2,3\:\:\mbox{mod $4$}$.
\item The $p$-torsion of $MString_*$ is trivial for every prime $p\ge 5$.
\item The $3$-torsion of $MString_*$ is annihilated by multiplication with $3$.
\end{enumerate}
The spin bordism groups $MSpin_*$ are calculated additively in \cite{MR0190939} (see Section \ref{udiqdqwdwqdqwdwqdwqd} for a table of $MSpin_*$).
We will use the following facts:
\begin{enumerate}
\item  $MSpin_k$ is a finite group for $k\equiv 1,2,3\:\:\mbox{mod $4$}$.
\item For all $k\ge 0$ the torsion $MSpin_{k,tors}\subseteq MSpin_k$ is a direct
sum of copies of $\Z/2\Z$.
\item $MSpin_{4m-1}=0$ for $m\le 9$.
\end{enumerate}
We get the following consequences for  the subgroup
$A_{4m-1}\subseteq MString_{4m-1}$:
\begin{enumerate}
\item For all $m\ge 1$, the group $A_{4m-1}\subseteq MString_{4m-1}$ contains 
 all $3$-torsion elements.
\item For all $m\ge 1$ we have $2 MString_{4m-1}\subseteq A_{4m-1}$. 
\item  We have $ A_{4m-1}=MString_{4m-1}$ for $m\le 9$.
\end{enumerate}

Because of our lack of   knowledge of $A_{4m-1}$, 
for an explicit description of $b^{top}$  Proposition \ref{udiqwdqwdqwdq} is a crucial observation.
It asserts that $b^{top}$ has a factorization
 $$\xymatrix{MString_{4m-1}\ar[r]^\sigma&tmf_{4m-1}\ar[d]^{b^{tmf}}&\\A_{4m-1}\ar@{_{(}->}[u]\ar[r]^{b^{top}}&T_{2m}} \ .
$$ 
In contrast to $MString$, the homotopy groups of $tmf$ are known   \cite{MR2508200}.  In Propositions \ref{uidodqwdqwd} and \ref{zwdqdwqd} we obtain a complete calculation of $b^{tmf}$. In particular, the $3$-torsion of $tmf_{4m-1}$ is detected completely by $b^{tmf}$.

This calculation implies that $b^{top}$ is non-trivial in arbitrarily high dimensions.
Let us explain some examples. 
\begin{itemize}
\item
There exists a $3$-torsion element
$x\in tmf_{27}$ which goes to the element with the name
$\nu\Delta\in tmf_{(3),27}$ under  localization at $3$. 
By Proposition \ref{zwdqdwqd} we know that $b^{tmf}(x)\not=0$.

By a result of Hopkins-Mahowald \cite[Theorem 6.25]{MR1989190} the $tmf$-valued Witten genus  $$\sigma:MString_*\to tmf_*$$ is surjective. We can choose an element $y\in A_{27}=MString_{27}$  such that $\sigma(y)=x$. 
 We then have 
$$ b^{top}(y)\stackrel{ \ref{udiqwdqwdqwdq}}{=}b^{tmf}(x)\not=0\ .$$
 
If $(M,\alpha)$ is a closed $27$-dimensional string manifold one can try to calculate the element
$\sigma([M,\alpha])\in  tmf_{(3),27}= (\Z/3\Z)\:\nu\Delta$. An answer in terms of characteristic numbers is given in Subsection \ref{ziddqwdqwdwq}.
Similar statements can be produced in higher dimensions using the $72$-periodicity of $tmf_{(3),*}$ given by the multiplication with  $\Delta^3$.
\item
  Let $x\in tmf_{192}$ be an element which goes to
$\Delta^8\in tmf_{(2),192}$ under localization at $2$. Then by Hopkins-Mahowald \cite[Theorem 6.25]{MR1989190} there exists a class $y\in MString_{192}$ such that $\sigma(y)=x$. Let $g\in A_3=MString_3$ be the generator (\ref{reuwrewrwer}).
 Note that $yg\in A_{195}$. The element
$\sigma( yg)\in tmf_{195}$ goes to $\nu \Delta^8\in tmf_{(2),195}$ under localization at $2$.
By Proposition \ref{uidodqwdqwd} the order of $b^{top}(yg)=b^{tmf}(\nu\Delta^8)$ is $8$.
\end{itemize}

\subsection{Open problems}

We close with stating some open problems:
\begin{enumerate}
\item What is $A_{4m-1}\subseteq MString_{4m-1}$? For which $m$ do we have an equality?
\item What is the image $\sigma(A_{4m-1})\subseteq tmf_{4m-1}$? For which $m$ do we have an equality? 
\item What is $b^{an}(x)\in T_{2m}$ for $x\in MString_{4m-1}\setminus A_{4m-1}$?
One could conjecture that
$$b^{an}=b^{tmf}\circ \sigma\ .$$

\end{enumerate}

\section{Three-manifolds}\label{gfdhfdf}

\subsection{A string bordism invariant in dimension $3$}\label{ziqdqwdwqd}

In the present subsection we give a geometric construction of a homomorphism
$$d:MString_{3}\to \Z/24\Z\ .$$
We will see in Corollary \ref{uqiduwqdqwdwqdd} that $d$ is an isomorphism.
 It is known that $MSpin_3=0$. Let $M$ be a closed spin three-manifold $M$.
Then we can  find a spin zero bordism
$Z$. We choose a connection $\nabla^{TM}$ on $TM$. It naturally induces a connection on the $Spin(3)$-principal bundle given by the spin structure.
We furthermore choose an extension of $\nabla^{TM}$ to a connection $\nabla^{TZ}$.

Let us fix a topological string structure $\alpha^{top}$ on $M$.
As a principal bundle can be equipped with a connection
a  topological string structure $\alpha^{top}$ on a spin bundle with connection can be refined to a geometric string structure $\alpha$. For details we refer to  \cite{waldorf}.
A geometric string structure $\alpha$ on $M$ gives rise  to a $3$-form $H_{\alpha}$
which satisfies $dH_\alpha=\frac{1}{2} p_1(\nabla^{TM})$  (this condition is non-vacuous in higher dimensions).
We can form the difference\footnote{This difference has  also been considered in the recent (and independent) paper
\cite{redden}, where it is put in relation with the $e$-invariant of Adams.}
$$d_Z(M,\alpha):=\frac{1}{2} \int_Z p_1(\nabla^{TZ})-\int_MH_{\alpha}\in \R\ .$$
\begin{lem}\label{zudqwdqwdqwd}
The real number $d_Z(M,\alpha)$ is independent of the choice of connections and the   geometric data of the string structure.
\end{lem}
\proof
First observe that the difference does not depend on the extension of the connection to the bordism $Z$. Indeed, given two extensions $\nabla^{TZ},\nabla^{TZ,\prime}$  we form the closed manifold
$W:=Z\cup_M -Z$ with the induced connection $\nabla^{TW}$ which restricts to
$\nabla^{TZ}$ and $\nabla^{TZ,\prime}$ on the two obvious copies of $Z$ in $W$.  
We then have
$$\int_Z p_1(\nabla^{TZ})-\int_Z p_1(\nabla^{TZ,\prime})=\int_W p_1(\nabla^{TW})=3\  \sign(W)\ .$$
But $\sign(W)=0$ because of the obvious orientation-reversing $\Z/2\Z$-symmetry of $W$.

Two connections on $TM$ can be joined by a connection $\nabla^{T(I\times M)}$ on the  cylinder $I\times M$.
Similarly,  two geometric string structures $\alpha$ and $\alpha^\prime$ refining the same underlying topological string structure can be connected by a geometric string structure $\tilde \alpha$ on $I\times M$.
Then we have by Stokes' theorem
\begin{equation}\label{zqwduqgwdqwdwqdwqd}
\frac{1}{2}\int_{I\times M}p_1(\nabla^{T(I\times M)})=\int_{I\times M} dH_{\tilde \alpha}=\int_M H_\alpha
-\int_M H_{\alpha^\prime}\ .\end{equation}
We choose $Z^\prime:=Z\cup_M (I\times M)$ as the zero bordism  for the primed choices. Then 
\begin{equation}\label{zqwduqgwdqwdwqdwqd1}\frac{1}{2}\int_{Z^\prime} p_1(\nabla^{TZ,\prime})-\frac{1}{2}\int_{Z} p_1(\nabla^{TZ})=\frac{1}{2}\int_{I\times M}p_1(\nabla^{T(I\times M)})\ .
\end{equation}
The assertion now follows from the combination of  (\ref{zqwduqgwdqwdwqdwqd}) and (\ref{zqwduqgwdqwdwqdwqd1}). 
Alternatively one could conclude the result from Lemma \ref{uefhiwefuwefw} below
and a continuity argument.
\hB

As a consequence of Lemma \ref{udidqwdqwdqwdqwdqwdqwd},  the integral $\int_{M} H_\alpha$ does not depend on the geometry of $\alpha$ so that we can also write it as 
$\int_{M} H_{\alpha^{\tiny top}}$. Moreover we can write $d_Z(M,\alpha)=d_Z(M,\alpha^{top})$.
The set of topological string structures on $M$ is a torsor under $H^3(M;\Z)\cong \Z$. We will write the action additively, see   (\ref{ufiwefwefwefw}).
Note that (see \cite{waldorf}) 
\begin{equation}\label{uiqdqwdwqdqwdwd54545}
\int_{M} H_{\alpha^{\tiny top}+x}=\int_M H_{\alpha^{\tiny top}}+\langle x,[M]\rangle\ ,\quad \forall x\in H^3(M,\Z)\ .
\end{equation}
Therefore
\begin{equation}\label{udiqdwqdqwd}
d_Z(M,\alpha^{top}+x)=d_Z(M,\alpha^{top})-\langle x,[M]\rangle\ .
\end{equation}

\begin{lem}\label{uefhiwefuwefw}
We have $d_Z(M,\alpha)\in \Z$.
\end{lem}
\proof
We let $\widehat{H\Z}^*$ denote the differential integral cohomology functor. 
It has first been introduced by  Cheeger-Simons \cite{MR827262}. For an axiomatic picture see \cite{bunke-2009}. Differential integral  cohomology is the home for geometric refinements  of Chern and Pontrjagin classes.
Let $$\frac{\hat p_1}{2}(\nabla^{TM})\in \widehat{H\Z}^4(M)$$ denote the 
lift introduced in \cite{MR827262} of the spin characteristic class $\frac{p_1}{2}(TM)$  (see  also \cite[Sec. 4.2]{bunke-20095}).  Its integral over $M$ is an element of $\R/\Z$. By the bordism formula for the evaluation of 
differential cohomology classes we have the equality 
$$[\frac{1}{2}\int_Z p_1(\nabla^{TZ})]=\int_M \frac{\hat{p}_1}{2}(\nabla^{TM})$$
in $\R/\Z$.
We also know (see e.g. \cite[Sec. 4.2]{bunke-20095}) that
$$a(H_\alpha)= \frac{\hat{p}_1}{2}(\nabla^{TM})\ ,$$
where $a:\Omega^3(M)/\im(d)\to \widehat{H\Z}^4(M)$ is one of the structure maps
of differential cohomology.
It follows that
$$[d_Z(M,\alpha^{top})]=[\frac{1}{2}\int_Z p_1(\nabla^{TZ})-\int_MH_{\alpha}]=0\in \R/\Z\ .$$
\hB

By the specialization of Hirzebruch's signature theorem to the four-dimensional case, the 
 signature of a closed oriented $4$-dimensional spin manifold $W$ is given by 
\begin{equation}\label{uidqwdqwdww}
\sign(W)=\frac{1}{3}\int_W p_1(\nabla^{TW})\ . 
\end{equation}
Similarly, by the specialization of the Atiyah-Singer index theorem   the index of the spin Dirac operator of a four-dimensional spin manifold $W$ is given by
\begin{equation}\label{uidqwdqwdww1}\ind(D_W)=-\frac{1}{24}\int_W p_1(\nabla^{TW})\ . \end{equation} 
In dimensions $8m+4$ the spin Dirac operator has an additional real symmetry which forces its index to be even.
It follows that
$$ \frac{1}{2}\int_W p_1(\nabla^{W})\equiv 0 \:\:\:\:\mbox{mod $24$}\ .$$

From the additivity of the signature
we conclude that the class
$$d(M,\alpha^{top}):=[d_Z(M,\alpha^{top})]\in \Z/24\Z$$
is independent of the choice of the zero bordism $Z$.

\begin{lem}\label{uiqodqwdqwdqwd}
$d(M,\alpha^{top})\in \Z/24\Z$ is a string bordism invariant.
\end{lem}
\proof
If $\alpha$ has an extension to a geometric string structure $\tilde \alpha$ over $Z$, then by Stokes' theorem
$$
\frac{1}{2}\int_Z p_1(\nabla^{TZ})=\int_Z dH_{\tilde \alpha}=
 \int_M H_\alpha\ .
$$
\hB 

It is now clear that $$d:MString_3\to \Z/24\Z\ ,\quad [ M,\alpha^{top}]\mapsto d(M,\alpha^{top})$$ is well-defined and a homomorphism.

\subsection{A generator for $MString_3$}

In this subsection we construct  a specific generator $g\in MString_3$
such that $d(g)=[1]\in \Z/24\Z$. It will be used in later calculations.
Let $S$ denote the sphere spectrum, and let  $\epsilon:S\to MString$ be the unit of the ring spectrum $MString$. It is known \cite[Thm 2.2.1]{MR1455523} that
$\epsilon:S_3\to MString_3$ is an isomorphism so that
$MString_3\cong S_3\cong \Z/24\Z$. 

We consider the sphere $S^3\subset \R^4$. As the boundary of the disc $D^4\subset \R^4$ it has a preferred orientation, spin structure and string structure $\alpha^{top}$.  We define
 \begin{equation}\label{reuwrewrwer}
g:=[S^3,\alpha^{top}-\ori_{S^3}] \in MString_3\ , 
\end{equation}
where $\ori_{S^3}\in H^3(S^3;\Z)$ is the orientation class of $S^3$.  
 It follows immediately from (\ref{udiqdwqdqwd}) that
$d(g)=[1]\in \Z/24\Z$. Since the order of $d(g)$ is $24$ we  see that $g\in MString_3$ is a generator. 
\begin{kor}\label{uqiduwqdqwdwqdd}
The homomorphism
$$d:MString_3\to \Z/24\Z\ ,\quad  [M,\alpha^{top}]\mapsto d(M,\alpha^{top})$$
is an isomorphism
\end{kor}
We will systematically extend this construction to higher dimensions in  Section \ref{uiudwdqwdqwdd}.

\subsection{Atiyah's canonical $2$-framing}\label{uidqwdqwdqw}

In this Subsection we discuss, as a side aspect of the main topic of the present paper,   an analogue for  spin manifolds  of an observation by Atiyah \cite{MR1046621} saying
that oriented three-manifolds have canonical $2$-framings.
We will show  that every closed  connected spin three-manifold has a canonical string structure.

Let us first explain the result of Atiyah \cite{MR1046621}. For all $n\ge 2$ the horizontal composition in the following diagram of Lie groups
$$\xymatrix{&&Spin(2n)\ar[d]\\SO(n)\ar@{.>}[urr]\ar[r]_/-.9em/{\diag}\ar[r]&SO(n)\times SO(n)\ar[r]&SO(2n)}
$$
has a unique lift indicated by the dotted arrow. This implies that
the  double $2V:=V\oplus V$ of  a $n$-dimensional real oriented vector bundle $V$ has a canonical spin structure. 

A $2$-framing of a closed oriented three-manifold $M$ is by definition a
spin-trivialization of  the double $2TM$ of its tangent bundle. 
Atiyah now considers an oriented  zero bordism $Z$ of $M$. 
The spin bundle $2TZ$ is trivialized  by the $2$-framing $\alpha$  at the boundary $\partial Z\cong M$. This trivialization refines the spin characteristic class
$\frac{p_1}{2}(2TZ)\in H^4(Z;\Z)$ to a relative cohomology class 
 $$\frac{p_1}{2}(2TZ,\alpha)\in H_c^4(Z,M;\Z)\ .$$
Atiyah  then observes that
$$\sigma(\alpha):=3\sign(Z)- \langle \frac{p_1}{2}(2TZ,\alpha),[Z,M]\rangle\in \Z$$ does not depend on the oriented zero bordism
$Z$. Furthermore, by changing $\alpha$, we can alter the right-hand side by any integer.
It follows that there is a unique $2$-framing $\alpha_0$ such that
$\sigma(\alpha_0)=0$. This is the canonical $2$-framing of $M$.

We now define the canonical topological string structure of a three-dimensional closed connected spin manifold $M$. Let $\alpha^{top}$ be any topological string structure on $M$. The combination
$$\sigma(M,\alpha^{top}):=3\sign(Z)-  2 d_Z(M,\alpha^{top})\in  \Z$$
is independent of the choice of the spin zero bordism $Z$  of $M$
(by the same argument as in \cite{MR1046621} which employs additivity of the signature and the formula (\ref{uidqwdqwdww})).
Its class
$$\sigma(M):=[\sigma(M,\alpha^{top})]=[3\sign(Z)-2d_Z(M,\alpha^{top})]=[\sign(Z)]\in  \Z/2\Z$$ is independent of the choice of the string structure $\alpha^{top}$ as well.  
We will see from the examples below that both possible values do occur.

\begin{prop}\label{uwifwefwef}
A closed  connected spin three-manifold $M$ has
a unique topological string structure $\alpha^{top}_0$ characterized by 
$\sigma(M,\alpha_0^{top})\in \{0,1\}$.
\end{prop}

\begin{ddd}\label{uiqdwqdqwdw}
The topological string structure uniquely characterized in Proposition \ref{uwifwefwef}
will be called the canonical string structure.
\end{ddd}

\subsection{Examples}\label{zuqdqwd}

In this subsection we calculate two examples which show that both possible  values in of the invariant $\sigma(M)\in \Z/2\Z$ and thus of $\sigma(M,\alpha^{top}_0)\in \{0,1\}$ of a connected three-dimensional spin manifold do really occur. 
\begin{enumerate}
\item  The three-dimensional compact group $SO(3)$ is diffeomorphic to $\R¸\P^3$. We thus  have
$$H^*(SO(3);\Z/2\Z)\cong \Z/2\Z[x]/(x^4)\ ,\quad |x|=1\ .$$ The manifold $SO(3)$ fits  into a principal bundle $SO(2)\to SO(3)\stackrel{\pi}{\to} S^2$. The circle $SO(2)$ acts by the adjoint representation on the quotient of Lie algebras $so(3)/so(2)$. We fix a basis element of the Lie algebra $so(2)$.
It determines a complex structure on $so(3)/so(2)$ by requiring that it rotates counterclockwise.
It furthermore  induces a fundamental vector field on $SO(3)$ which trivializes the vertical
 bundle $T^v\pi$. The vertical bundle therefore gets a spin and a topological string structure $\alpha_{T^v\pi}^{top}$.  This spin structure induces the bounding spin structure on all fibers  of $\pi$.

Furthermore, using the identification $TS^2\cong SO(3)\times_{SO(2)}so(3)/so(2)$, we fix a complex structure on $TS^2$. This complex structure determines the orientation of $S^2$. 
Since $H^1(S^2;\Z/2\Z)=0$ and $w_2(TS^2)=0$ it can be refined to a spin structure in a unique way. Furthermore, there exists a unique
topological string structure $\alpha_{S^2}^{top}$ on the spin manifold $S^2$.
Up to homotopy there is a unique  splitting 
\begin{equation}\label{udiqdwqdwqd}
TSO(3)\cong T^v\pi\oplus \pi^*TS^2 .
\end{equation}
The spin structures of the summands induce a spin structure $s_0$ on $SO(3)$.
Since $H^{1}(SO(3);\Z/2\Z)\cong \Z/2\Z$ there is actually a second spin structure $s_1$.
Since the restriction of $H^{1}(SO(3);\Z/2\Z)\to H^{1}(\pi^{-1}(x);\Z/2\Z)$ is an isomorphism for all $x\in S^2$  the spin structure $s_1$ induces the non-bounding spin structure on every fiber of $\pi$.
We continue with the spin structure $s_0$ and choose the topological string structure $\alpha^{top}$ such that (\ref{udiqdwqdwqd}) becomes an isomorphism
of string bundles.

The trivialization of $T^v\pi$ induces a geometric string structure $\alpha_{T^v\pi}$ refining $\alpha^{top}_{T^v\pi}$ with $H_{\alpha_{T^v\pi}}=0$. We further choose a geometric string structure $\alpha_{S^2}$ which refines $\alpha^{top}_{S^2}$.
We have $H_{\alpha_{S^2}}=0$ since $S^2$ is two-dimensional. We choose the geometric string structure  $\alpha$ refining $\alpha^{top}$ as  the sum of $\alpha_{T^v\pi}$ and $\pi^*\alpha_{S^2}$. Then we have $H_\alpha=0$.

The action of $SO(2)$ on $so(3)/so(2)$ fixes a metric and  induces an action on the unit disc $D^2$.
If we fix the point $1\in D^2$, then we can identify $SO(2)$ with the boundary of $D^2$, the orbit of $1$. In this way we get an identification  $SO(3)\cong \partial Z$ with the boundary of the four-dimensional manifold
 $$Z:=SO(3)\times_{SO(2)} D^2\ .$$
The $D^2$-bundle $q:Z\to S^2$ exhibits this manifold as a fiberwise zero bordism of $SO(3)$.
 It can be identified with the unit-disc bundle in the tangent bundle $TS^2$.
The vertical bundle of $q$ is therefore given by $T^vq\cong q^*TS^2$. Up to homotopy
we have a unique decomposition
$$TZ\cong T^vq\oplus q^*TS^2\cong q^*(TS^2\oplus TS^2)\ .$$  
It gives $TZ$ a complex structure and therefore an orientation. Since the projection
$q:Z\to S^2$ is a homotopy equivalence and $H^1(S^2;\Z/2\Z)=0$ we see that the spin structure on $TZ$ induced by the spin structure on $TS^2$ and this decomposition is the unique one. The restriction of the spin structure on $T^vq$ to each fiber of $\pi$
is obviously the bounding one. It follows that the restriction of the spin structure on $Z$ to its boundary is the spin structure $s_0$ fixed above. Hence $Z$ is a spin zero bordism of  $SO(3)$.

We have identified $SO(3)$ with the unit sphere bundle in $TS^2$. This unit sphere bundle
is the same as the orthonormal complex frame bundle. Therefore the complex bundle
$\pi^*TS^2\cong T^vq_{|SO(3)}$ is canonically trivialized.

We can choose a complex connection $\nabla^{T^vq}$ which is compatible
with this trivialization, and we define
$\nabla^{TZ}$ as the sum of $\nabla^{T^vq}$ and $q^*\nabla^{TS^2}$, where
$\nabla^{TS^2}$ is any complex connection of $TS^2$.
Then we have
$$p_1(\nabla^{TZ})=-c_2(\nabla^{TZ})=-c_1(\nabla^{T^vq})\wedge  q^*c_1(\nabla^{TS^2})\ .$$
By a standard calculation\footnote{Without any explicit calculation of integrals this can be seen as follows.
We consider the complex  line bundle $D^2\times \C\to  D^2$ with the trivialization
over $S^1$ given by the tautological section $s:S^1\to \C$ on the boundary. Using this
trivialization we glue it with the trivial bundle $\bar D^2\times \C$, where $\bar D$  has the opposite complex structure,  in order to get a holomorphic bundle $L\to S^2\cong \C\P^1$. We choose a connection on
$D^2\times \C$ which is compatible with the trivialization at the boundary. This connection can then be extended to $\nabla^L$ as a trivial connection across $\bar D$. 
The section $z:D^2\to D^2\times \C$ extends by $1$ to a holomorphic section of $L$. 
Furthermore, the constant section  $1:D^2\to  D^2\times \C$ can be extended by the holomorphic function $\bar z$ to $\bar D$.
It is easy to see that these two generate the two-dimensional space of holomorphic  sections
of $L\to \C\P^2$. By the Riemann-Roch theorem
$2=1+\int_{\C\P^1} c_1(L)=1+\int_{D^2} c_1(\nabla^{L})$ and therefore
$\int_{D^2} c_1(\nabla^{L})=1$.
} we get
$\int_{Z/S^2} c_1(\nabla^{T^vq})=1$.
It follows that
$$\int_{Z}p_1(\nabla^{TZ})=-\int_{S^2}c_1(\nabla^{TS^2})=-2\ .$$
This gives $d_Z(SO(3),\alpha^{top})=-1$.

Using again that  $q:Z\to S^2$ is a homotopy equivalence we
see that $H_2(Z;\Z)\cong \Z$ is generated by the orientation of $S^2$.
The self intersection of this class, which geometrically can be  considered as the zero section of $TS^2$, is the Euler characteristic  $\chi(S^2)=2$. It is positive so that
 $\sign(Z)=1$.
We conclude that
$$\sigma(SO(3),\alpha^{top})=3+2=5\ .$$
It follows that $$\sigma(SO(3))=1\ .$$
The canonical string structure of $SO(3)$ with the spin structure $s_0$ fixed above  is given by  $$\alpha_{0}^{top}=\alpha^{top}-2\ \ori_{SO(3)}\ .$$
 \item
\newcommand{\T}{\mathbb{T}}
We now consider the torus $\T^3\cong \R^3/\Z^3$.
This representation as a quotient fixes a trivialization of $T\T^3$
and therefore a spin and geometric string  structure $\alpha$ with $H_\alpha=0$.
The induced spin structure on each circle subgroup of $\T^3$ is bounding.
We write $\T^3\cong S^1\times \T^2$ and can consider $S^1$ as the spin boundary of the disc $D^2$. Therefore 
$Z:=D^2\times \T^2$ is a spin zero bordism of $\T^3$. 
We choose $\nabla^{TZ}$ such that respects this decomposition and
the factor $\T^2$ is flat. Then $p_1(\nabla^{TZ})=0$.
It follows that $d_Z(\T^3,\alpha^{top})=0$. The intersection form 
on  $H_2(Z;\Z)$ vanishes so that $\sign(Z)=0$. We conclude that
$$\sigma(\T^3)=0$$
for the spin structure chosen above,  and that the canonical string structure of $\T^3$ with this spin structure is given by
$$\alpha_0^{top}=\alpha^{top}\ .$$
\end{enumerate}

 \section{Invariants from the Witten genus}\label{uiudwdqwdqwdd}

\subsection{Introduction}

In Section \ref{gfdhfdf}, using  geometric string structures, we have defined a homomorphism $$d:MString_3\to \Z/24\Z$$ which turned out to  be an isomorphism. 
In the present section, with the construction of the homomorphisms $b^{geom}:A_{4m-1}\to T_{2m}$, we generalize this   to higher dimensions.

It is interesting  to observe that using the index theorem of Atiyah-Patodi-Singer \cite{MR0397797}  we can give an alternative expression for  $d$ which does not involve the zero bordism $Z$. This parallels the treatment of the $e$-invariant of Adams given in  \cite{MR0397797}.
We will actually work in the slightly different setting of index theory for manifolds with boundary which has been developed in \cite{MR2191484}, and which we will    briefly review in Subsection \ref{eoiwewe}.

We consider a closed string three-manifold $(M,\alpha^{top})$ and a spin zero bordism $Z$.
Let us assume that $\nabla^{TM}$ is the Levi-Civita connection associated to a Riemannian metric on $M$, and that $\nabla^{TZ}$ is the Levi-Civita connection for an extension of that metric to $Z$ with a product structure.
These geometric structures turn the manifolds $M$ and $Z$ into  geometric manifolds $\cM$ and $\cZ$. If we choose a (real) taming $\cM_t$, then we get a boundary tamed manifold $\cZ_{bt}$.
The index $\ind(\cZ_{bt})$ of the Fredholm operator associated to a boundary tamed manifold can be calculated by the formula \cite[Thm. 2.2.18]{MR2191484}.
It follows from the presence of a real structure that $\ind(\cZ_{bt})$ is even. In the case at hand we have
$$\ind(\cZ_{bt})=-\frac{1}{24}\int_Z p_1(\nabla^{TZ}) + \eta (\cM_t)\ ,$$
where $\eta(\cM_t)$ is the eta-invariant of the  tamed manifold $\cM_t$, see Subsection \ref{eoiwewe}.
 Therefore we have the following equality in $\R/\Z$:
$$[\frac{1}{2}\int_Z p_1(\nabla^{TZ})]=[12\eta(\cM_t)]\ .$$
We furthermore choose a geometric refinement $\alpha$ of the topological string structure  $\alpha^{top}$
based on the spin connection induced by $\nabla^{TM}$.
Then we get the expression
$$d(M,\alpha)=[12\eta (\cM_t)-\int_M H_\alpha]\in \Z/24\Z$$
which is now intrinsic to $M$.
With the construction of $b^{an}:MString_{4m-1}\to T_{2m}$ (see Definition \ref{uiqdwqdqwdwqd} ) we generalize this 
  analytic formula  to higher dimensions.

\subsection{Tamings and $\eta$-invariants}\label{eoiwewe}

In this subsection we recall the necessary language of local index theory for manifolds with boundary. The main reference for the set up  is  \cite{MR2191484} which specialized to manifolds with boundary can be considered as a variant of \cite{MR0397797}. 
If $M$ is a $4m-1$-dimensional closed Riemannian spin manifold and $V\to M$ is a real vector bundle with metric $h^V$  and connection $\nabla^V$, then we can form the twisted Dirac operator $D_M\otimes V$ which acts on sections of the bundle $S(M)\otimes_\R V$, where $S(M)\to M$ is the spinor bundle.   A taming of $D_M\otimes V$ is by definition a selfadjoint operator $Q$ acting on sections of $S(M)\otimes V$ and given by a smooth integral kernel such that $D_M\otimes V+Q$ is invertible. If $m$ is odd, then the spinor bundle has a quaternionic
structure. If the taming respects this quaternionic structure, then we call it a real taming. If $m$ is even, then $S(M)$ has a real structure, and a real taming should respect this structure. The obstruction against the existence of a real taming is
$\ind(D_M\otimes V)\in KO_{4m-1}$. Since $KO_{4m-1}=0$ for all $m\ge 1$  such real tamings always exist. Following the notation introduced in \cite{MR2191484} we will call a Riemannian spin manifold a geometric manifold $\cM$, and we denote such a manifold with an additional choice of a geometric bundle $\bV:=(V,h^V,\nabla^V)$ by $\cM\otimes \bV$. Finally,  this data together with  a choice of a real taming will be denoted by $(\cM\otimes \bV)_t$ and called a tamed manifold. 
The $\eta$-invariant
$\eta((\cM\otimes \bV)_t)\in \R$
of a tamed manifold is defined
by
$$\eta((\cM\otimes \bV)_t):=\frac{-1}{\sqrt{\pi}}\int_0^\infty \Tr (D_M+Q)e^{-t^2(D_M+Q)^2}dt\ .$$
These definitions can be extended to graded or virtual bundles $V$ in a natural way.

\subsection{The Witten genus}\label{udqwidqwdqwdqwd444}

In this subsection we recall the Witten genus $R:MSpin\to KO[[q]]$ 
and introduce several related formal power serieses.
We use the correspondence between an even formal power series
$\phi\in A[[x^2]]$ and a power series
$K_\phi\in A[[p_1,p_2,\dots]]$ given by
$$\prod_{i=1}^\infty \phi(x_i)=K_\phi(p_1,p_2,\dots)\ , \quad
\sum_{i\ge 0} p_i=\prod_{i=1}^\infty (1+x_i^2)\ ,$$
where  $A$ is some commutative  $\Q$-algebra, e.g. $A=\Q[[q]]$ in the example which follows. 
We let (see \cite[Ch. 6.3]{MR1189136})
$$\phi_W(x,q)\in \Q[[q]][[x]]$$ be given by
\begin{eqnarray}
\phi_W(x,q)&:=&\frac{\frac{x}{2}}{\sinh(\frac{x}{2})} \prod_{n=1}^\infty\frac{(1-q^n)^2}{(1-q^ne^x)(1-q^ne^{-x})}\nonumber\\
&=&\exp\left[\sum_{k=2}^\infty \frac{2}{(2k)!}G_{2k} x^{2k}\right] e^{G_2(q)x^2}\ .\label{ujkwefef6e7}
\end{eqnarray}
Here 
$$G_{2k}:=-\frac{B_{2k}}{4k}+\sum_{n=1}^\infty \sigma_{2k-1}(n)q^n\in \Q[[q]]$$
with the Bernoulli numbers $B_{2k}$ and the sum $\sigma_{2k-1}(n)$ of the $(2k-1)$'th-powers of the positive divisors of $n$. 
For $k\ge 2$ the power series 
$G_{2k}$ is the $q$-expansion of a modular form of weight $2k$.  We have $G_{2k}\in \cM_{2k}^\Q$, since the constant term of the formal power series $G_{2k}$ is  rational, while all higher terms are integral.

We will also need the series
$$\theta_W(x,q):=\exp\left[\sum_{k=2}^\infty \frac{2}{(2k)!}G_{2k} x^{2k}\right]\in  \Q[[q]][[x]]\ .$$
We define 
\begin{equation}\label{udiqwdqwdqwd3}
\Phi:=K_{\phi_W}\in \Q[[q]][[p_1,p_2,\dots]]\ ,\quad \Theta:=K_{\theta_W}\in  \Q[[q]][[p_1,p_2,\dots]]\ .
\end{equation}

 For $k\ge 1$ we let
 $$N_{2k}(p_1,\dots)=\sum_{j=1}^\infty x_j^{2k}\in \Q[[p_1,\dots]]$$
be the Newton polynomials.
Then we can write
$$\Phi(p_1,\dots)=\exp\left[\sum_{k=2}^\infty \frac{2}{(2k)!}G_{2k} N_{2k}(p_1,\dots )\right]e^{G_2p_1}\in \Q[[q]][[p_1,\dots]]\ .$$
We define
$$\tilde \Phi(p_1,\dots,):=\exp\left[\sum_{k=2}^\infty \frac{2}{(2k)!}G_{2k} N_{2k}(p_1,\dots )\right]   \sum_{j=1}^\infty   \frac{G_2^j p_1^{j-1}}{j!}\in \Q[[q]][[p_1,\dots]] \ .$$
More systematically
\begin{equation}\label{udiqdqwdqwdqwd}
\tilde \Phi =\Theta \frac{e^{G_2p_1}-1}{p_1}\ .
\end{equation}

We now turn to the definition of the Witten genus.
For a real $l$-dimensional vector bundle $V\to B$ we define the formal power series of virtual bundles
\begin{equation}\label{zfquefhefwefwef}
R(V):=\sum_{n\ge 0} q^n R_n(V)= \prod_{k\ge 1} (1-q^k)^{2l}\bigotimes_{k\ge 1} Sym_{q^k}(V\otimes_\R \C)\ ,
\end{equation} 
where
$$Sym_q(V):=\bigoplus_{n\ge 0} q^n Sym^n(V)$$
is the generating power series of the symmetric powers $Sym^n(V)$ of $V$.
For a manifold $Z$ we then have 
\begin{equation}\label{zduzwtduqwdwd}\hA(TZ)\cup \ch(R(TZ))  = \Phi(p_1(TZ),\dots)\ .\end{equation}
If $\nabla^{TZ}$ is a connection on $TZ$, then we will write
$$\Phi(\nabla^{TZ}):=\Phi(p_1(\nabla^{TZ}),\dots)\in \Omega^{4*}(Z)[[q]]\ ,$$
and we will use a similar convention for $\tilde \Phi$ and $\Theta$.

The homomorphism 
$$R:MSpin_{4m}\to \Z[[q]]\cong KO[[q]]_{4m}$$
is given by
\begin{equation}
\label{udiqdqwdqwdqwdw}R([M])=\kappa_m \ind(D_M\otimes R(TM))=\kappa_m \langle \Phi(p_1(TM),p_2(TM),\dots),[M]\rangle 
\end{equation}
(see \cite{MR1189136}).
In the case of odd $m$ the factor $\kappa_m=\frac{1}{2}$ takes care of the factor $2$ on the lower horizontal arrow of the commutative diagram
$$\xymatrix{KO_{4m}\ar[d]^{\cong}\ar[r]&K_{4m}\ar[d]^\cong\\
\Z\ar[r]^2&\Z}\ ,$$
since we take the index of the twisted Dirac operator as a complex operator.

The above construction can be refined to   a multiplicative natural transformation of homology theories $$R:MSpin\to KO[[q]]\ .$$
Let $f:M\to X$ be a map from a closed $l$-dimensional spin manifold $M$ to a space $X$ representing the  bordism class $[M,f]\in MSpin_l(X)$.
If we choose a Riemannian metric on $M$, and $V\to M$ is a real vector bundle with connection, then we can define
the real $K$-homology class
$[D_M\otimes V]\in KO_{l}(M)$ of the Dirac operator of $M$ twisted by $V$.
It is independent of the choice of the geometric structures.
The value of $$R:MSpin_l(X)\to KO[[q]]_l(X)$$ on $[M,f]\in MSpin_l(X)$ is then given by 
$$R([M,f]):=\sum_{n\ge 0} q^n f_*[D_M\otimes R_n(TM)]\in KO_l(X)[[q]]=KO[[q]]_l(X)\ .$$

The formal power series defined in (\ref{zfquefhefwefwef}) is multiplicative in the sense that $R(V\oplus W)=R(V)\otimes R(W)$. This easily implies that the transformation $R$
is multiplicative.

\subsection{A string bordism invariant in dimension $4m-1$}\label{dzuidwqdqwd}

Let $m\ge 1$ and consider a $4m-1$-dimensional closed string manifold $(M,\alpha^{top})$. We choose a Riemannian metric on $M$ and get a geometric manifold $\cM$. The Riemannian metric gives a Levi-Civita connection $\nabla^{TM}$
which together with the trivial connection on the trivial bundle $\Theta_\R=M\times \R$ induces a connection on the virtual  bundles $R_n(TM\oplus \Theta_\R)$ for all $n\ge 0$.
We choose real tamings $(\cM\otimes R_n(TM\oplus \Theta_\R))_t$ for all $n\ge 0$.
Finally we choose a geometric refinement $\alpha$ of the topological string structure $\alpha^{top}$. It gives rise to a form $H_\alpha\in \Omega^3(M)$ satisfying
$dH_\alpha=\frac{1}{2}p_1(\nabla^{TM})$.

\begin{ddd}\label{uiqdwqdqwdwqd}
We define the formal power series 
\begin{equation}\label{formal}\tilde b^{an}(\cM,\alpha,t):= 2\kappa_m\int_M H_{\alpha}\wedge \tilde \Phi(\nabla^{TM})+\kappa_m\sum_{n\ge 0} q^n \eta((\cM\otimes R_n(TM\oplus \Theta_\R))_t) \in \R[[q]] \ ,\end{equation}
where 
$$\kappa_m=\left\{\begin{array}{cc}
1&m\equiv 0(2)\\
 \frac{1}{2}&m\equiv 1(2)
 \end{array}\right.\ ,$$
and  we choose tamings which are compatible with the real structure.
\end{ddd}
The last entry in the list of variables of $\tilde b^{an}$  indicates the dependence of the formal power series on the choice of the tamings.

Recall the Definition \ref{uzifefwfwefewfwef} of the group $T_{2m}$.  

 \begin{ddd}\label{uzduiwqdqwd}
We let   $$ b^{an}(M,\alpha^{top}):=[\tilde b^{an}(\cM,\alpha,t)]\in T_{2m}$$ denote the class in $T_{2m}$ represented by  the formal power series $\tilde b^{an}(\cM,\alpha,t)$. 
\end{ddd}
This notation is justified by the following Lemma.

\begin{lem}
 The class $b^{an}(M,\alpha^{top})\in T_{2m}$ is an invariant of the string bordism class of the $4m-1$-dimensional string manifold $(M,\alpha^{top})$. In particular, it is independent of the choice of geometric structures and the taming involved in the definition of $\tilde b^{an}(\cM,\alpha,t)$.
\end{lem}
\proof
Let $(Z,\tilde \alpha)$ be a string bordism between
$(M,\alpha)$ and $(M^\prime,\alpha^\prime)$. 
Then we choose a Riemannian metric on $Z$ which extends the given metrics on $M$ and $M^\prime$ with product structures. We have decompositions of geometric bundles $TZ_{|M}\cong TM\oplus \Theta_{\R}$ and $TZ_{|M^\prime}\cong TM^\prime\oplus \Theta_{\R}$, where the trivial summands correspond to the normal bundle.
Therefore,
for all $n\ge 0$ the geometric bundles
$R_n(TZ)$ extend the geometric bundles $R_n(TM\oplus \Theta_\R)$ and $R_n(TM^\prime\oplus \Theta_\R)$  The tamings $(\cM\otimes R_n(TM\oplus \Theta_\R))_t$ and $(\cM^\prime\otimes R_n(TM^\prime\oplus \Theta_\R))_{t^\prime}$ induce a boundary taming
$(\cZ\otimes R_n(TZ))_{bt}$.
The index theorem \cite[Thm. 2.2.18]{MR2191484} for boundary tamed manifolds gives
\begin{eqnarray*}\lefteqn{
\ind(\cZ_{bt}\otimes R_n(TZ))}&&\\&&=\int_{Z} \hA(\nabla^{TZ})\wedge \ch(\nabla^{R_n(TZ)}) + \eta((\cM^\prime\otimes R_n(TM^\prime\oplus \Theta_\R))_{t^\prime})-\eta((\cM\otimes R_n(TM\oplus \Theta_\R))_t)\ .
\end{eqnarray*}
If $m$ is odd, then $\ind(\cZ_{bt})$ is even. Using that (see (\ref{zduzwtduqwdwd})), 
\begin{equation}\label{uidqwdqwdqwd}
\sum_{n\ge 0} q^n  \hA(\nabla^{TZ})\wedge \ch(\nabla^{R_n(TZ)})=\Phi(\nabla^{TZ})\ .
\end{equation}
we get the following equality in $T_{2m}$.
 \begin{eqnarray}\lefteqn{\hspace{-2cm}
\left[\kappa_m\sum_{n\ge 0}q^n\big(\eta((\cM^\prime\otimes R_n(TM^\prime\oplus \R))_{t^\prime})-\eta((\cM\otimes R_n(TM\oplus \R))_t)\big)\right] }&& \nonumber\\&=&
[-\kappa_m\int_Z \Phi(\nabla^{TZ})]  \hspace{4cm}\label{klwefwef}\end{eqnarray}
Using Stoke's theorem, (\ref{udiqdqwdqwdqwd}), and $d H_{\tilde \alpha}=p_1(\nabla^{TZ})$ we calculate
\begin{eqnarray}
2\int_{M^\prime} H_{\alpha^\prime}\wedge \tilde \Phi(\nabla^{TM^\prime}) -
2\int_{M } H_{\alpha }\wedge \tilde \Phi(\nabla^{TM})&=&\int_Z 2dH_{\tilde \alpha} \wedge  
\tilde \Phi(\nabla^{TZ})\nonumber\\
&=&\int_Z \Theta(\nabla^{TZ}) e^{G_2 p_1(\nabla^{TZ})}-\int_Z\Theta(\nabla^{TZ})\nonumber\\
&=&\int_Z \Phi(\nabla^{TZ})-\int_Z\Theta(\nabla^{TZ})\label{uihzdqwieerqweqwe}\ .
\end{eqnarray}
Note that $\int_Z\Theta(\nabla^{TZ})\in \cM^\R_{2m}$ so that
\begin{equation}\label{uwdqwdqwd}\left[2\kappa_m\int_{M^\prime} H_{\alpha^\prime}\wedge \tilde \Phi(\nabla^{TM^\prime}) -
2\kappa_m\int_{M } H_{\alpha }\wedge \tilde \Phi(\nabla^{TM})\right]=
\left[\kappa_m\int_Z \Phi(\nabla^{TZ})\right] 
\end{equation}
in $T_{2m}$.
If we combine
(\ref{uwdqwdqwd}) and (\ref{udiqdqwdqwdqwd}), then  we get
$b^{an}(M,\alpha)=b^{an}(M^\prime,\alpha^\prime)$. 

The independence of $b^{an}(M,\alpha^{top})$ from the geometric structures and tamings 
follows from the bordism invariance since we can connect two choices of geometric structures and tamings by corresponding structures on
 a cylinder $Z=[0,1]\times M$ which provides a bordism.
\hB 

The invariant $b^{an}(M,\alpha^{top})$ is clearly additive under disjoint unions of string manifolds.
We have therefore defined a homomorphism
\begin{equation}\label{udiudwqdqwdqwdwqd}
b^{an}:MString_{4m-1}\to T_{2m}\ ,\quad [M,\alpha^{top}]\mapsto b^{an}(M,\alpha^{top})\ .
\end{equation}
The main goal of the present paper is to understand the homotopy theoretic meaning
of this global analytic construction.

\subsection{A geometric expression}\label{ioqdqwd}

Since the construction of $b^{an}$ involves spectral invariants of Dirac operators,
a direct evaluation of $b^{an}$   is complicated. Our analysis of $b^{an}$ relies on the comparison with a homotopy theoretic version $b^{top}:A_{4m-1}\to T_{2m}$. The comparison between
$b^{an}_{|A_{4m-1}}$ and $b^{top}$ is achieved  via an intermediate construction
$$b^{geom}:A_{4m-1}\to T_{2m}$$
using differential geometry which we will present in the present subsection.

 We consider a closed string manifold $(M,\alpha^{top})$ of dimension $4m-1$ which represents an element in $[M,\alpha^{top}]\in A_{4m-1}\subseteq MString_{4m-1}$.  Then we can choose a spin zero bordism  $Z$ of $M$. We choose a connection $\nabla^{TZ}$ on $TZ$ and let $\nabla^{TM}$ be its restriction to $M\cong \partial Z$.
Furthermore, we let $\alpha$ be a geometric string structure which  refines
$\alpha^{top}$  based on the spin connection on $TM$ induced by $\nabla^{TM}$.
Then we can define the formal power series
 \begin{equation}\label{w8e9fwefwfwefwf}
\tilde b^{geom}(M,\alpha,\nabla^{TZ}):=2\kappa_m \int_M H_\alpha\wedge \tilde \Phi(\nabla^{TM})-\kappa_m \int_Z \Phi(\nabla^{TZ}) \in \R[[q]]\ .\end{equation}

\begin{lem}\label{udidqwdqwd}
In $T_{2m}$ we have the equality
$$[\tilde b^{geom}(M,\alpha,\nabla^{TZ})]=b^{an}(M,\alpha^{top})\ .$$
\end{lem}
\proof
We first observe that $[\tilde b^{geom}(M,\alpha,\nabla^{TZ})]$ does not depend  on the choice of the
connection  $\nabla^{TZ}$ and the corresponding geometric refinement $\alpha$ of $\alpha^{top}$. 
Let $\nabla^{TZ}_i$ and $\alpha_i$, $i=0,1$,  be  two     choices of geometric structures. Then we can find a connection $\nabla^{TW}$  on $W:=[0,1]\times Z$
which restricts to $\nabla^{TZ}_i$ on $\{i\}\times Z$. Furthermore we can find a string structure $\tilde \alpha$ on $V:=[0,1]\times M$ which connects the string structures $\alpha_i$, $i=0,1$.  From
(\ref{uihzdqwieerqweqwe})
we get
$$
2\int_M H_{\alpha_1}\wedge \tilde \Phi(\nabla^{TM}_1)- 2 \int_M H_{\alpha_0}\wedge \tilde \Phi(\nabla_0^{TM}) = \int_{V} \Phi(\nabla^{TV}) -\int_{V}\Theta (\nabla^{TV})\ .
$$
The last term belongs to $\cM_{2m}^\R$ so that
 \begin{equation}\label{d77uqdidqwdwqd}
\left[2\kappa_m \int_M H_{\alpha_1}\wedge \tilde \Phi(\nabla^{TM}_1)-2\kappa_m \int_M H_{\alpha_0}\wedge \tilde \Phi(\nabla_0^{TM})\right]=\left[\kappa_m\int_{V} \Phi(\nabla^{TV})\right]
\end{equation}
in $T_{2m}$. 
Again by  Stoke's theorem we have
\begin{eqnarray}
0&=&\int_{W} d\Phi(\nabla^{TW})\nonumber
=\int_{\partial W} \Phi(\nabla^{TW})\\
&=&\int_{Z}\Phi(\nabla^{TZ}_1)-\int_{Z}\Phi(\nabla^{TZ}_0)-\int_{V}\Phi(\nabla^{TV}_0)\ .\label{zgsuzqs}
\end{eqnarray}
If we combine (\ref{zgsuzqs}) and  (\ref{d77uqdidqwdwqd}), then we get
$[\tilde b^{geom}(M,\alpha_1,\nabla^{TZ}_1)]=[\tilde b^{geom}(M,\alpha_0,\nabla^{TZ}_0)]$.

We now choose a Riemannian metric on $Z$ with product structure near the boundary which extends the metric of $M$. Then we get the geometric manifold $\cZ$. The taming $(\cM\otimes R_n(TM\oplus \Theta_\R))_t$ induces a boundary taming  $(\cZ\otimes R_n(TZ))_{bt}$, and by the index formula 
$$\kappa_m\int_{Z}\hA(\nabla^{TZ})\wedge \ch(\nabla^{R_n(TZ)})+\kappa_m\eta((\cM\otimes R_n(TM\oplus \R))_t)=\kappa_m\ind((\cZ\otimes R_n(TZ))_{bt})\in \Z$$
for all $n\ge 0$.
Using (\ref{uidqwdqwdqwd}) we get the equality of classes
$$\left[-\kappa_m \int_Z \Phi(\nabla^{TZ})\right]=\left[\kappa_m\sum_{n\ge 0}q^n \eta((\cM\otimes R_n(TM\oplus \R))_t)\right]$$
in $T_{2m}$. Combining this with (\ref{formal}) and (\ref{w8e9fwefwfwefwf})
we conclude that 
$$[\tilde b^{geom}(M,\alpha,\nabla^{TZ})]=[\tilde b^{an}(\cM,\alpha,t)]=b^{an}(M,\alpha^{top})\ .$$
\hB 
As a consequence of Lemma \ref{udidqwdqwd} we have constructed a homomorphism
$$b^{geom}:A_{4m-1}\to T_{2m}$$ which maps $[M,\alpha^{top}]\in A_{4m-1}$ to the class in $T_{2m}$ represented by the formal power series
$b^{geom}(M,\alpha,\nabla^{TZ})$ for some choices of the spin zero bordism $Z$ and the geometric structures. Of course, we have  the equality 
$$b_{|A_{4m-1}}^{an}=b^{geom}\ .$$

\subsection{A topological expression}\label{udqwidqwdqwdqwd}

In this subsection we retain the assumptions of Subsection \ref{ioqdqwd}
 and transform the geometric expression $[\tilde b(M,\alpha,\nabla^{TZ})]=b^{geom}(M,\alpha^{top})$   into a purely topological
formula. Let $[0,1)\times M\hookrightarrow
Z$ be a collar and consider a function
$\chi\in C^\infty_c[0,1)$ with $\chi\equiv 1$ near $0$.
Using the normal variable of the collar this function can be transported to the collar. Its extension by zero gives a smooth function on $Z$ which we will also denote by $\chi\in C^\infty(Z)$.

Let $p:[0,1)\times M\to M$ denote the projection.
The form  $d(\chi p^* H_\alpha)$
can be considered as a closed form on $Z$ supported near $M$.
We now define the closed form
$$\tilde p_1(\nabla^{TZ},\alpha):=p_1(\nabla^{TZ})-2d(\chi p^*H_\alpha)\in \Omega^4_c(Z)\ .$$
It represents a relative cohomology class \begin{equation}\label{uidquwdqwdqwdqwdqwdqwdq}
\tilde p_1(TZ,\alpha^{top}):=[\tilde p_1(\nabla^{TZ},\alpha)]\in H^4(Z,M;\R)\ .
\end{equation}
 
We define 
$$\hat  b^{geom}(M,\alpha^{top},Z):=\langle -\kappa_m \tilde \Phi(TZ)\cup \tilde p_1(TZ,\alpha^{top}), [Z,M]\rangle\in \R[[q]]\ .$$

\begin{lem}\label{zdqwud}
In $T_{2m}$ we have the equality
$$b^{geom}(M,\alpha^{top})=[\hat b^{geom}(M,\alpha^{top},Z)]\ .$$
\end{lem}
\proof
We write
\begin{eqnarray*}
\Theta(\nabla^{TZ}) &=&\Phi_W(\nabla^{TZ}) -\tilde \Phi_W(\nabla^{TZ})   \wedge p_1(\nabla^{TZ}) \\&=&\Phi_W(\nabla^{TZ}) -\tilde \Phi_W(\nabla^{TZ})  \wedge \tilde p_1 -2\tilde \Phi_W(\nabla^{TZ})   \wedge d(\chi H_\alpha)\\&=&\Phi_W(\nabla^{TZ}) -\tilde \Phi_W(\nabla^{TZ}) \wedge  \tilde p_1 -2d (\tilde \Phi_W(\nabla^{TZ})  \wedge \chi H_\alpha)\ .
\end{eqnarray*}
Since $\int_Z \Theta(\nabla^{TZ})\in \cM_{2m}^\R$  and
$$\int_{Z} 2d (\tilde \Phi (\nabla^{TZ}) \wedge \chi H_\alpha)=2\int_M \tilde \Phi(\nabla^{TZ}) \wedge H_\alpha$$
 we conclude that
$$[b^{geom}(M,\alpha,\nabla^{TZ})]=[- \kappa_m\int_Z   \tilde \Phi(\nabla^{TZ})\wedge  \tilde p_1]=[\hat b^{geom}(M,\alpha,Z)]$$
in $T_{2m}$.

\subsection{Calculation in the case $m=1$}

In the present Subsection we explain the relation between the three-dimensional constructions in Section \ref{gfdhfdf} with the higher-dimensional theory of the present section.  
We have $\kappa_1=\frac{1}{2}$. In (\ref{reuwrewrwer}) we have fixed a generator
$g\in MString_3$ in the form
$$g=[S^3,\alpha^{top}-\ori_{S^3}]=[S^3,\alpha^{top}-\ori_{S^3}]-[S^3,\alpha^{top}]\ ,$$
where $\alpha^{top}$ is the topological string structure coming from the representation of $S^3$ as the boundary of the disc $D^4\subset \R^4$ with its natural topological string structure.

 We use formula (\ref{w8e9fwefwfwefwf}) in order to calculate $b^{geom}(g)\in T_2$.
Because of the structure of the representative of $g$ as a difference, the contribution of the zero bordism $Z$ drops out. In the three-dimensional case only the 
cohomological-degree $0$ part $\tilde \Phi_{[0]}=G_2$ of $\tilde \Phi$ contributes.  Using (\ref{uiqdqwdwqdqwdwd54545}) we get 
$$b^{geom}(g)=b^{geom}([S^3,\alpha^{top}-\ori_{S^3}])-b^{geom}([S^3,\alpha^{top}]) =\left[-2\kappa_mG_2\langle \ori_{S^3},[S^3]\rangle\right] =\left[-G_2\right]  \ .$$
Note that
$$G_2=-\frac{1}{24}+ q+\dots\ .$$
The higher terms of the $q$-expansion are all integral.
This implies 
\begin{equation}\label{dzuwdzquwd}
b^{geom}(g)= [\frac{1}{24}]\in T_2=\R/\Z[[q]]\ .
\end{equation}
Therefore  $b^{geom}(g)$ has order $24$,  and we see that
$b^{geom}:MString_3\to T_2$ is injective.
This shows in particular that the homomorphism $b^{geom}=b^{an}$  
is non-trivial at least in some examples.

Let $j:\Z/24\Z\to \R/\Z$ be the inclusion such that
$j([1]):=[\frac{1}{24}]$. Then we get the following relation between the three-dimensional construction $d$ as in Corollary \ref{uqiduwqdqwdwqdd}
and $b^{geom}$:
$$\xymatrix{MString_3\ar[r]^d\ar[d]^{b^{geom}}&\Z/24\ar[d]^j\Z\\T_2\ar[r]^{p_0}&\R/\Z}\ ,
$$
where $p_0:T_2=\R/\Z[[q]]\to \R/\Z$ takes the constant coefficient.

\section{The secondary  index theorem}

\subsection{Construction of $b^{top}$}\label{qdwqdwqd}

In this subsection  for all $m\ge 1$ we give a  homotopy theoretic construction of a homomorphism 
$$b^{top}:A_{4m-1}\to T_{2m}\ .$$

We let $ko\to KO$ denote a connective cover of the real $K$-theory spectrum.
Then we get an induced connective cover 
 $$c:ko[[q]]\to KO[[q]]\ .$$

We let $tmf$ denote the connective spectrum of topological modular forms
constructed by Goerss, Hopkins and Miller \cite{MR1989190}, \cite{MR2125040}.  It fits into the following commutative diagram
of multiplicative transformations
$$\xymatrix{&&ko[[q]]\ar[d]^c\\MString\ar[d]^j \ar[r]^\sigma&tmf\ar[r]^W\ar@{.>}[ur]^w&KO[[q]]\\MSpin\ar@{..>}@/^0.5cm/[rruu]^r\ar[rru]_{R}&}\ .$$
 
The map $\sigma$ will be called the $tmf$-valued Witten genus. The $tmf$-valued Witten genus  and  the factorization
 $W\circ \sigma=R\circ j:MString \to KO[[q]]$ (for $R$ see Subsection \ref{udqwidqwdqwdqwd444}) have  been constructed by Ando, Hopkins, Rezk  \cite{ahr}. Since $tmf$ is connective we have a unique factorization
of $W$ over a ring map $w:tmf\to ko[[q]]$.
Similarly, since $MSpin$ is connective, we have a unique ring map $r$ which lifts $R$.

For an abelian  group $G$ we let $MG$ denote the  associated Moore spectrum. For a spectrum $X$
we write $XG:=X\wedge MG$. The inclusion $\Z\to \Q$ induces a map $M\Z\to M\Q$, and we get an exact triangle
$$ M\Z\to M\Q \to M\Q/\Z\to \Sigma M\Z\  $$in the stable homotopy category. 
By smashing with $X$ it induces an exact triangle 
$$X\to X\Q \to X\Q/\Z\to \Sigma X$$
which is 
functorial in $X$.
For a spectrum $X$ we write $X^{l}:=\Sigma^l X$ for all $l\in \Z$.
We define the spectra $F$ and $G$ to fit into exact triangles
$$G^{-1}\to tmf\stackrel{w}{\to} ko[[q]]\to G$$
and
$$F^{-1}\stackrel{\nu}{\to} MString \stackrel{j}{\to} MSpin\to F\ .$$
We choose a map $\bar\sigma:G\to F$ such that the following diagram becomes a morphism of distinguished triangles
$$\xymatrix{F^{-1}\ar[r]^(.4){\nu}\ar[d]^{-\Sigma^{-1}(\bar\sigma)}&MString\ar[d]^\sigma\ar[r]^j&MSpin\ar[d]^r\ar[r]&F\ar[d]^{\bar\sigma}\\ G^{-1}\ar[r]&tmf\ar[r]^w&ko[[q]]\ar[r]&G}\ .$$

We now consider the following commutative  diagram:
\begin{equation}\label{udqidqwdqwdqwd}
\xymatrix{&&&MSpin\Q^{-1}\ar[rr]^{-\Sigma^{-1}(r_\Q)}\ar[dd]&&ko[[q]]\Q^{-1}\ar[dd]&\\
&&MSpin^{-1}\ar[ur]\ar[rr]^/-2em/{-\Sigma^{-1}(r)}\ar[dd]&&ko[[q]]^{-1}\ar[dd]\ar[ur]&&\\
&&&F\Q^{-1}\ar[rr]^/-2em/{-\Sigma^{-1}(\bar \sigma_\Q)}\ar[dd]&&G\Q^{-1}\ar[dd]&&\\
&&F^{-1}\ar[rr]_/-2em/{-\Sigma^{-1}(\bar\sigma)}\ar[ur]\ar[dd]\ar@/^0.0cm/@{..>}[rrru]^a&&G^{-1}\ar[dd]\ar[ur]_b&&\\
&&&MString\Q\ar[rr]^/-2em/{\sigma_\Q}\ar[dd]&&tmf\Q\ar[dd]&\\
S^{4m-1}\ar@{.>}@/_2.2cm/[ddrr]^0\ar[ddrr]\ar@{~>}@/^3cm/[uuuuurrr]_\omega\ar@{-->}@/_0.3cm/[uuurrr]_/1em/{F_{Z_\Q}}\ar@{-->}@/^0.3cm/[urrr]_/+1.6em/{(M,\alpha^{top})_\Q}\ar@{..>}[rr]_{(M,\alpha^{top})}\ar@{..>}[uurr]^{F_Z}&&MString\ar[ur]\ar[rr]^/-2em/{\sigma}\ar[dd]^j&&tmf\ar[ru]\ar[dd]^/-2em/w&&\\
&\ar@{::>}[dl]_Z&&MSpin\Q\ar[rr]^/-2em/{r_\Q}&&ko[[q]]\Q&\\
&&MSpin\ar[ur]\ar[rr]^r&&ko[[q]]\ar[ru]&&}\ .\end{equation}
The maps $b, \sigma_\Q$ and $\bar \sigma_\Q$ are the smash products of $\id_ {G^{-1}},\sigma $ and $ \bar \sigma$ with the canonical map $S=M\Z\to M\Q$, respectively. The map $a$ is defined as the composition $a:=-b\circ \Sigma^{-1}(\bar \sigma)$.   
Note that $a$ does not depend on the choice of $\bar \sigma$ because by construction
$$\xymatrix{F^{-1}\ar[r]^/-0.5em/\nu\ar[d]^{a}&MString\ar[d]\ar[r]^j&MSpin\ar[d]^r\ar[r]&F\ar[d]^{-\Sigma(a)}\\ G\Q^{-1}\ar[r]&tmf\Q\ar[r]^{w_\Q}&ko[[q]]\Q\ar[r]&G\Q}$$
is a morphism of triangles.
It determines $a$ up to elements coming from
$[MString,ko[[q]]]\Q]$. Since $MString$ is rationally even we have  $[MString,ko[[q]]\Q^{-1}]=0$.

We calculate the homotopy groups of $G\Q$ from the exact sequence
$$\dots \to tmf\Q_j \to ko[[q]]\Q_j \to G\Q_{j}\to tmf\Q_{j+1} \to ko[[q]]\Q_{j+1} \to \dots\ .$$
In particular, since $tmf$ is rationally even, we get the short exact sequences 
$$0\to tmf\Q_{4m}\to ko[[q]]\Q_{4m}  \to G\Q_{4m}\to  0\ .$$ 
For $m\ge 0$ we have $ko[[q]]_{4m}\cong \Z[[q]]$, so that  $ko[[q]]\Q_{4m}\cong \Z[[q]]\otimes\Q\subset \Q[[q]]$. Note that
$$\cM^\Q_{2m}\subset \Z[[q]]\otimes\Q\subset \Q[[q]]\ ,$$
and by \cite{MR1989190} the subgroup $\cM^\Q_{2m}$ coincides with  the image of 
$$\sigma_\Q:tmf\Q_{4m}\to ko[[q]]\Q_{4m}\ .$$ We therefore
get an identification
$$G\Q_{4m}\cong \frac{\Z[[q]]\otimes \Q}{\cM^\Q_{2m}}\ .$$

Let $(M,\alpha^{top})$ be a closed string manifold of dimension $4m-1$.
The associated string bordism class $[M,\alpha^{top}]\in MString_{4m-1}$ is represented by a  map
$(M,\alpha^{top}):S^{4m-1}\to MString$ as indicated in diagram (\ref{udqidqwdqwdqwd}). 
We assume that $$[M,\alpha]\in A_{4m-1}\subseteq MString_{4m-1}\ .$$
Then we can  choose a spin zero bordism
  $Z$ of $M$.  
This zero bordism can be interpreted as a zero homotopy $Z$ 
of the composition $$S^{4m-1}\stackrel{(M,\alpha^{top})}{\longrightarrow} MString\stackrel{j}{\to} MSpin$$ indicated in (\ref{udqidqwdqwdqwd}). It thus gives rise  to a map
$F_Z:S^{4m-1}\to F^{-1}$ representing the
element $[F_Z]\in F_{4m}$.
We have
$$a([F_Z])\in G\Q_{4m}\cong \frac{\Z[[q]]\otimes \Q}{\cM^\Q_{2m}}\ .$$
Note that
$[F_Z]$ is determined by $[M,\alpha^{top}]\in MString_{4m-1}$ up to elements in the image of $MSpin_{4m}\to F_{4m}$.
By a diagram chase we see that these go into the image of
$ko[[q]]_{4m}\to G\Q_{4m}$, i.e. of
$$\Z[[q]]\to  \frac{\Z[[q]]\otimes \Q}{\cM^\Q_{2m}}\ .$$
The natural chain of inclusions
$$\Z[[q]]\otimes \Q\to \Q[[q]]\to \R[[q]]$$
induces an inclusion
$$\frac{\Z[[q]]\otimes \Q}{\Z[[q]]+\cM^\Q_{2m}}\to \frac{\R[[q]]}{\Z[[q]]+\cM_{2m}}=T_{2m}\ .$$
In the following we will use these inclusions implicitly.
By construction,  the class $[a([F_Z])]\in T_{2m}$ only depends on the class $[M,\alpha^{top}]\in A_{4m-1}$. Moreover, it depends on this bordism class additively thus justifying the following definition.
\begin{ddd}\label{uzidqwdqwdwqdwd}
We define the homomorphism
\begin{equation}\label{uidqwdqwd}b^{top}:A_{4m-1}\to T_{2m}\ ,\quad [M,\alpha^{top}]\mapsto[a([F_Z])] \ .\end{equation}
\end{ddd}

\subsection{The index theorem}\label{uidqwdqwdqwdqwdee}

In this subsection we show the secondary index theorem $b^{top}=b^{geom}$.
It is very similar to the secondary index theorem $e^{an}=e^{top}$ for the $e$-invariant stated in  \cite[Introduction]{bunke-2008}. Let us explain the analogy in greater detail.
The primary invariant
for a map $S^k\to S$  in the case of the $e$-invariant was induced by the unit $x:S\to MU$. The $e$-invariant was obtained from a zero homotopy of the composition $S^{2m-1}\to S\to MU$, seen through the eyes of $K$-theory.
The construction can be visualized in the diagram
$$\xymatrix{&\overline{MU}^{-1}\ar[d]\ar[r]&K\Q/\Z^{-1}\\S^{2m-1}\ar@{-->}[urr]_/2em/{e(x)}\ar@{.>}[ur]\ar[dr]^0\ar[r]^x&S\ar[d]&\\
&MU&}\ .$$
Here $\overline{MU}$ is the cofiber of the unit, the dotted arrow is induced by the zero homotopy, and the map $\overline{MU}\to K\Q/\Z$
is explained in \cite{bunke-2008}.
In the present case the primary invariant for a map $x:S^k\to MString$ is induced by  the map $j:MString\to MSpin$.
The secondary invariant $b^{top}$ measures a zero homotopy of the composition
$S^{4m-1}\to MString\to MSpin$  (which exists since we assume that $x\in A_{4m-1}$) by means of the Witten genus $r$ with values in $ko[[q]]$.
The visualization is 
$$\xymatrix{&F^{-1}\ar[d]\ar[r]^a&G\Q^{-1}\\S^{4m-1}\ar@{-->}[urr]_/2em/{\hat b^{top}(x)}\ar@{.>}[ur]^{F_Z}\ar[dr]^0\ar[r]^x&MString\ar[d]&\\
&MSpin&}\ ,$$
where the class
$$\hat b^{top}(x)=a\circ  F_Z\in G\Q_{4m}\cong \frac{\Z[[q]]\otimes \Q}{\cM^\Q_{2m}}$$ represents $b^{top}(x)\in T_{2m}$.
The reason that we can not go directly to $ko[[q]]\Q/\Z$   is that
in contrast to $S_{2m}\otimes\Q=0$ for $m\ge 1$ in the case of the $e$-invariant, in our situation  we have  $MString_{4m}\otimes \Q\not=0$ in general.

\begin{theorem}\label{e2e7uiwqdqwdwqd}
For all $m\ge 1$ we have the following equality of homomorphisms
$$b^{geom}=b^{top}:A_{4m-1}\to T_{2m}\ .$$
\end{theorem}
\proof
We let $(M,\alpha^{top}):S^{4m-1}\to MString$  be a map 
which represents a class $[M,\alpha^{top}]\in A_{4m-1}$.
It  produces a  closed $4m-1$-dimensional string manifold, also denoted by $(M,\alpha^{top})$, via the Thom-Pontrjagin construction and Lemma \ref{udiqdqwdwqd}. 
 Since
$MString$ is rationally even the element
$[F_{Z_\Q}]\in F\Q_{4m}$ maps  to zero in $MString\Q_{4m}$.
This implies the existence of a lift $\omega:S^{4m-1}\to MSpin\Q^{-1}$ as indicated in (\ref{udqidqwdqwdqwd}).
By a diagram chase we see that the formal power series
$-\Sigma^{-1}(r_\Q)(\omega)\in ko[[q]]\Q_{4m}\cong \Z[[q]]\otimes \Q$ represents 
$b^{top}([M,\alpha^{top}])$.

Let us give the geometric interpretation of $\omega$.
Since $MString\Q_{4m-1}=0$   there exists a $4m$-dimensional string manifold $(V,\beta^{top})$ and a natural number $L\in \nat$ such that
$L (M,\alpha^{top})=\partial (V,\beta^{top})$, where by $L (M,\alpha^{top})$ we denote the disjoint union of $L$ copies of $(M,\alpha^{top})$.
We define a spin manifold $LZ\cup_M -V$, where
$-V$ denotes the spin manifold $V$ with the opposite orientation and spin structure. Then
$\omega=\frac{1}{L}[L Z\cup_M -V]$. 

We choose Riemannian metrics on $Z$ and $V$ with product structures which glue nicely.
Then we have associated Levi-Civita connections on the tangent bundles of the manifolds. 
By (\ref{uiqdqwdqwdwqdwqd}) or (\ref{udiqdqwdqwdqwdw})  we have
$$-\Sigma^{-1}(r_\Q)(\omega)=\frac{\kappa_m}{L}\int_{LZ\cup_M -V} \Phi(\nabla^{T(LZ\cup_M -V)})\in \Z[[q]]\otimes\Q\cong k[[q]]\Q_{4m}\ .$$
Using a geometric refinement $\beta$ of the  string structure $ \beta^{top}$
with restriction $\alpha$ to $M$
we can write by (\ref{udiqdqwdqwdqwd})
$$ \Phi(\nabla^{TV})=2d H_\beta\wedge \tilde \Phi(\nabla^{TV})+
  \Theta(\nabla^{TV})\ .$$
By Stokes' theorem we get
\begin{eqnarray*}-\Sigma^{-1}(
r_\Q)(\omega)&=&\frac{\kappa_m}{L}\int_{LZ\cup_M -V} \Phi(\nabla^{T(LZ\cup_M -V)})\\&=&\kappa_m\int_Z \Phi(\nabla^{TZ})-2\kappa_m\int_MH_\alpha\wedge  \tilde \Phi(\nabla^{TM})-\frac{\kappa_m}{L} \int_V\Theta(\nabla^{TV})\ .
\end{eqnarray*}
The last integral belongs to $\cM^\R_{2m}$, and the first two term together equal $\tilde b^{geom}(M,\alpha,\nabla^{TZ})$, compare with (\ref{w8e9fwefwfwefwf}). 
This shows that 
$b^{top}([M,\alpha])=b^{geom}([M,\alpha])$.
\hB 

\subsection{Factorization over $tmf$}\label{uiffwefwefewf}

In this subsection we show that 
$$b^{top}:A_{4m-1}\to T_{2m}$$ defined in Definition \ref{uzidqwdqwdwqdwd} admits  a factorization over a homomorphism
$$b^{tmf}:tmf_{4m-1 } \to T_{2m}\ .$$
For the construction of $b^{tmf}$ we consider the following diagram of  vertical and horizontal fiber sequences
\begin{equation}\label{diagr2}
\xymatrix{&&&tmf\Q^{-1}\ar[d]\\
&&\:\:\:\:\:ko[[q]]^{-1}\ar[d]\ar[r]&\:\:{}^{\hat x_\Q}\:\:ko[[q]]\Q^{-1}\ar[d]\\&
\:\:\:\: G\Q/\Z^{-2}\ar[r]\ar[d]&\:\:{}^{\tilde x}\:\:G^{-1}\ar[r]\ar[d]&\:\:{}^{\tilde x_\Q}\:\:G\Q^{-1} \ar[d]\\&
\:\:{}^{x_{\Q/\Z}}\:\: tmf\Q/\Z^{-1}\ar[d]^{\bar w}\ar[r]&\:\:{}^x\:\:tmf\ar[d]^w\ar[r]&\:\:\:{}^{x_\Q}\:\: tmf\Q\ar[d]\\\:\:{}^z\:\: ko[[q]]\Q^{-1}\ar[r]&\:\:{}^{\bar w(x_{\Q/\Z})}\:\: ko[[q]]\Q/\Z^{-1}\ar[r]&\:\:\:\:\: ko[[q]]\ar[r]&\:\:\:\:\:ko[[q]]\Q}\ .\end{equation}
The script-size symbols denote elements which will be chased through the diagram during the following discussion.
Let $x\in tmf_{4m-1}$. We are going to construct $b^{tmf}(x)\in T_{2m}$ by a diagram chase.
Since $ko[[q]]_{4m-1}=0$ we can choose a lift
$\tilde x\in G_{4m}$. Its image $\tilde x_\Q\in G\Q_{4m}$ maps to $x_\Q=0$ in $tmf\Q_{4m-1}=0$. Therefore we can further choose a lift $\hat x_\Q\in ko[[q]]\Q_{4m}\cong \Z[[q]]\otimes \Q$ which represents
by definition $b^{tmf}(x)\in T_{2m}$. 

Let us check that
$b^{tmf}(x)$ is well-defined independently of the choice of the lifts. The lift
$\tilde x$  is defined uniquely up to elements
in the image $ko[[q]]_{4m}\to G_{4m}$. These are mapped to integral power serieses in
$ko[[q]]\Q_{4m}\cong \Z[[q]]\otimes \Q$. The lift $\hat x_\Q$ is well-defined up to elements in the image of $tmf\Q_{4m}\to ko[[q]]\Q_{4m}$, i.e. the space of power series $\cM_{2m}^\Q$.
Since in $T_{2m}$ integral and modular power serieses are factored out,
it follows that $b^{tmf}(x)$ is independent of the choices.

\begin{prop}\label{udiqwdqwdqwdq}
For every $m\ge 1$
we have a factorization 
 $$\xymatrix{MString_{4m-1}\ar[r]^/0.5em/\sigma&tmf_{4m-1}\ar[d]^{b^{tmf}}&\\A_{4m-1}\ar@{_{(}->}[u]\ar[r]^{b^{top}}&T_{2m}}\ .
$$ 
 \end{prop}
\proof
We chase in the diagrams (\ref{udqidqwdqwdqwd}) and (\ref{diagr2}).
Let $y\in A_{4m-1}\subseteq MString_{4m-1}$ and $x=\sigma(y)\in tmf_{4m-1}$.
Then we choose a lift $y_F\in F_{4m}$. The element $\tilde x:=-\Sigma^{-1}(\bar \sigma)(y_F)\in G_{4m}$ can serve as a lift of $x$ in the construction of $b^{tmf}$. The image
$\tilde y_\Q\in F\Q_{4m}$ can further be lifted to
$\hat y_\Q\in MSpin\Q^{-1}_{4m}$.  On the one-hand the element $-\Sigma^{-1}(r_\Q)(\hat y_\Q)\in ko[[q]]\Q_{4m}\cong \Z[[q]]\otimes \Q$
by construction represents $b^{top}(y)$. On the other hand, it can serve
as a lift $\hat x_\Q$ in the construction of $b^{tmf}(x)$ and thus
represents $b^{tmf}(x)$. 
\hB 

Using the diagram (\ref{diagr2}), we can give the following alternative description of $b^{tmf}(x)$ which will be employed in the Adams spectral sequence calculations in the Subsections \ref{uwiefwefwefwefewfwfef} and \ref{qwddwqdqwdoiuqwdoiwqd} below.  Since
$x\in tmf_{4m-1}$ maps to $x_\Q=0$ we can choose a lift $x_{\Q/\Z}\in tmf\Q/\Z_{4m}$. 
We let $\bar w(x_{\Q/\Z})\in ko[[q]]\Q/\Z_{4m}$ be its image under the map $\bar w:=w\wedge \id_{M\Q/\Z}:tmf\Q/\Z\to ko[[q]]\Q/\Z$ induced by $w $. Note that there is a natural projection map
$p:ko[[q]]\Q/\Z_{4m}\cong \Z[[q]]\otimes \Q/\Z \to T_{2m}$.
\begin{lem}\label{duqwiduqwdqwdqwd}
We have the equality
$$-p(\bar w({x_{\Q/\Z}}))=b^{tmf}(x)\ .$$
\end{lem}
\proof
We can  lift $\bar w(x_{\Q/\Z})$ 
to an element $z\in ko[[q]]\Q_{4m}$. 
By a diagram chase we see that
$$z+\hat x_\Q\in \im(tmf\Q_{4m}\to ko[[q]]\Q_{4m})+\im(ko[[q]]_{4m}\to ko[[q]]\Q_{4k})
$$
so that $\hat x_\Q$ and $-z$ represent the same element in $T_{2m}$. \hB

\section{Some details about Adams spectral sequences}

In this section we collect some facts and small results about 
 the generalized Adams spectral sequence in connection with
$\Q/\Z$-versions of spectra. This material will be used in the analysis of the invariant $b^{tmf}:tmf_{4m-1}\to T_{2m}$ in later sections.

\subsection{The $E_1$-term and the cobar complex}

Let us recall some details of the construction of the generalized Adams spectral sequence
${}^Y E_*^{*,*}(X)$ for a spectrum $X$ and a commutative ring spectrum $Y$ such that $Y_*Y$ is a flat $Y_*$-module.
As a first step we form the  Adams resolution of the sphere spectrum
$S$ \cite[Def. 2.2.10]{ravenel}
$$\xymatrix{S\ar[d]&\bar Y^{-1}\ar[d]\ar[l]&\bar Y^{-1}\wedge \bar Y^{-1}\ar[l]\ar[d]&\bar Y^{-1}\wedge \bar Y^{-1}\wedge \bar Y^{-1}\ar[l]\ar[d]&\dots \ar[l]\\Y\ar@{.>}[ur]& Y\wedge \bar Y^{-1}\ar@{.>}[ur]& Y\wedge\bar Y^{-1}\wedge\bar Y^{-1}\ar@{.>}[ur]&Y\wedge \bar Y^{-1}\wedge \bar Y^{-1}\wedge \bar Y^{-1}\ar@{.>}[ur]&}\ .$$
 Then we obtain the   standard $Y$-based Adams resolution of $X$   by smashing this resolution of $S$ with  $X$.

Note that $Y_*X$ is a comodule over the Hopf algebroid $(Y_*,Y_*Y)$. 
It gives rise to the reduced cobar complex
$\bar C^*(Y_*X)$, see \cite[A.1.2.11]{ravenel}. There is a decomposition $\overline{Y_*Y}\oplus Y_*\cong Y_*Y$ of $Y_*$-modules (actually there are two such decompositions
induced by the left or right unit of the Hopf algebroid), where $\overline{Y_*Y}\subset  Y_*Y$ is the kernel of the counit. Hence  the reduced cobar complex is a subcomplex of the cobar complex
$C^*(Y_{*}X)\cong Y_*Y^{\otimes_{Y_*} *}\otimes_{Y_*} Y_*X$.
 According to \cite[Prop. 2.11]{ravenel} the $E_1$-term of the generalized Adams spectral sequence is given by 
$$({}^Y E_1^{*,*}(X),d_1)\cong (\bar C^*(Y_{*}X),d^{cobar})\ .$$

\subsection{Application to $tmf$}\label{zudwzdqwd}

In this subsection we recall the relevant example of a certain ring spectrum $Y$
introduced by Ravenel \cite{MR737778}.
The $K$-theory functor $K^0$ is represented by the $H$-space $\Z\times BU$. We have an equivalence
$\Omega (\Z\times BU)\cong \Omega BU\cong U$.
This $H$-space represents the functor $K^{-1}$. Hence
$K^{-2}$ is represented by the $H$-space $\Omega U$.
The Bott periodicity transformation $K^{-2}\to K^0$ is represented by a map of $H$-spaces
$\Omega U\to \Z\times BU$. By restriction we get a 
map of $H$-spaces 
\begin{equation}\label{mmm1}\Omega U(4)\to \Z\times BU\end{equation}

To a map $\xi:X\to \Z\times BU$ we can associate a Thom spectrum $X^\xi$. If $\xi$ is a map of $H$-spaces, then $X^\xi$ is a ring spectrum which is commutative if $\xi$ deloops twice. 

It is known that the map (\ref{mmm1}) is a two-fold loop map and we 
let $Y$ be the associated commutative ring spectrum. 

The spectrum $tmf$ has the characterizing property that the cobar complex
$C^*(Y_{*}tmf)$ is the cobar complex of the Weierstrass Hopf algebroid \cite[14.5]{rezk}. 
In particular every group in the complex  $C^*(Y_{*}tmf)$ is  torsion-free.
The same is true for the reduced cobar complex, and  thus for the $E_1$-terms ${}^Y E_1^{*,*}(tmf_{(p)})$ of the $p$-localization $tmf_{(p)}$ of $tmf$ for all primes $p$.

\subsection{$\Q/\Z$-theory}

 For abelian groups $A,B$ we denote
$A*B:=Tor_1^\Z(A,B)$. One can identify
\begin{equation}\label{uidqwdqwdqwd55345345}A*\Q/\Z\cong A_{tors}\subseteq A\end{equation} 
in a canonical way, where
$A_{tors}\subseteq A$ denotes the torsion subgroup of $A$. 
For a spectrum $X$ we  have a functorial short exact sequence
\begin{equation}\label{ttrrz}0\to X_*\otimes \Q/\Z \to X\Q/\Z_*\to X_{*-1}* \Q/\Z \to 0\ .\end{equation}
The $Y$-based Adams resolution of $X\Q/\Z$ is obtained from the $Y$-based Adams resolution
of $X$ by smashing with the Moore spectrum $M\Q/\Z$. 
Therefore we get short exact sequences
$$0\to {}^Y E^{s,t}_1(X)\otimes \Q/\Z\to {}^Y E^{s,t}_1(X\Q/\Z)\to {}^Y E_1^{s,t-1}(X)*\Q/\Z\to 0\ .$$

The differentials of the $E_1$-term of the $Y$-based generalized Adams spectral sequence are induced by maps of spectra  defined before smashing with $M\Q/\Z$. Hence, from the functoriality of  (\ref{ttrrz}) in $X$,  we get a short exact sequences of complexes whose cohomology groups form long exact sequences involving  constituents of the $E_2$-terms of the generalized Adams spectral sequences of $X$ and $X\Q/\Z$.

 Now we assume that ${}^Y E_1^{*,*}(X)$ is torsion-free. Then
we have an isomorphism of complexes  $${}^Y E^{*,*}_1(X)\otimes \Q/\Z\cong {}^Y E^{*,*}_1(X\Q/\Z)\ .$$
>From the universal coefficient formula we deduce  exact sequences
$$0\to {}^Y E_2^{s,t}(X)\otimes \Q/\Z\to {}^Y E_2^{s,t}(X\Q/\Z)\stackrel{\partial}{\to} {}^Y E_2^{s+1,t}(X)* \Q/\Z\to 0\ .$$
Now we  consider the distinguished triangle in the stable homotopy category
$$  X\to X\Q\to X\Q/\Z\stackrel{h}{\to} \Sigma X  \ .$$
If $Y_*(h)=0$, then by  the geometric boundary theorem \cite[Thm. 2.3.4]{ravenel} for $r\ge 2$ there exist connecting maps
$$\delta_r:{}^Y E^{s,t}_r(X\Q/\Z)\to {}^Y E_r^{s+1,t}(X)$$
which are compatible with the differentials of the spectral sequences. Furthermore they are   filtered versions of $$h_*:X\Q/\Z_{*+1}\to X_*\ .$$

Assume again that ${}^Y E_1^{s,t}(X)$ is torsion free for all $s,t$. The case $s=0$ 
yields  $Y_*(h)=0$ so that the geometric boundary theorem applies. 
 In this way, using (\ref{uidqwdqwdqwd55345345}), we can compare the two maps
$$\xymatrix{{}^Y E_2^{s,t}(X\Q/\Z)\ar[r]^{\partial}\ar@/_0.6cm/[rr]_{\delta_2}&{}^Y E_2^{s+1,t}(X)_{tors}\ar[r]^i& {}^Y E_2^{s+1,t}(X)}\ ,$$
where $i$ denotes the canonical embedding.
\begin{lem}
We have $i\circ \partial=\delta_2$.
\end{lem}
\proof
Let $(A,\Gamma)$ be a Hopf algebroid   such that $\Gamma$ is a flat $A$-module, and let 
\begin{equation}\label{uzdiwdqwdqwd}
0\to M\to N\to Q\to 0
\end{equation} be an exact sequence of $(A,\Gamma)$-comodules. Then the associated complex of reduced cobar complexes $$0\to \bar C^*(M)\to \bar C^*(N)\to \bar C^*(Q)\to 0$$ is exact, too.  The boundary operator $$\Ext^s_{(A,\Gamma)}(A,Q)\to \Ext^{s+1}_{(A,\Gamma)}(A,M)$$ in the associated long exact sequence 
 is the $\Ext_{(A,\Gamma)}$-multiplication
with the extension class of the sequence (\ref{uzdiwdqwdqwd}) in
$\Ext^1_{(A,\Gamma)}(Q,M)$.
We now apply this to the exact complex
\begin{equation}\label{udiqwdqwd}
0\to Y_*X\to Y_*X\Q\to Y_*X\Q/\Z\to 0
\end{equation}
of $Y_*Y$-comodules.
The Adams $E_1$-terms give the exact sequences of reduced cobar complexes.
The map $i\circ \partial$ is the boundary operator of this sequence. The map $\delta_2$ is by
\cite[Thm 2.3.4 (a)]{ravenel} the $\Ext^1_{(A,\Gamma)}(Q,M)$-multiplication with the extension class of (\ref{udiqwdqwd}).
\hB

\subsection{Application to $tmf$}\label{udiqwdwqd}

In this short subsection we specialize the facts obtained above to the case of $tmf$.
We let $Y:=\Omega U(4)^\xi$ be as in Subsection \ref{zudwzdqwd}. Since
${}^Y E_1^{s,t}(tmf_{(p)})$ is torsion-free for all $s,t$ we get a short exact sequence
$$0\to {}^Y E^{s,t}_2(tmf_{(p)})\otimes \Q/\Z\to {}^Y E_2^{s,t}(tmf_{(p)}\Q/\Z)\stackrel{\partial}{\to} {}^Y E_2^{s+1,t}(tmf_{(p)})_{tors}\to 0\ .$$
Since $\partial$ is induced by $\delta_2$ it is compatible with  the boundary operators of the spectral sequences and a filtered version of
$$tmf_{(p)}\Q/\Z_{t-s}\to   tmf_{(p),t-s-1}\ .$$

\subsection{Application to $ko$}\label{uiqdqwdqwdqwd234324234}

As before we let $ko$ denote the connective real $K$-theory spectrum.
Replacing $U(4)$ by $U(2)$ in the construction of the spectrum $Y$ in Subsection  \ref{zudwzdqwd}
one gets a commutative ring spectrum $Y^\prime$. The Hopf algebroid $(Y^\prime_*,Y^\prime_*Y^\prime)$ is  flat. It can be used to calculate
$ko_*$. The Hopf algebroid 
\begin{equation}\label{ifowfwef}
(Y^\prime_*ko,Y^\prime_*Y^\prime\otimes_{Y^\prime_*}Y^\prime_*ko)\cong (\Z[b,c],\:\Z[b,c,r])
\end{equation}
with right unit
$$\left(\begin{array}{c}b\\c
\end{array}\right) \mapsto \left(\begin{array}{c}b+2r\\c+br+r^2 
\end{array}\right)$$
has been identified in \cite{hill}. Furthermore, in this reference it has been shown that the Hopf algebroid (\ref{ifowfwef}) is equivalent to the Hopf algebroid
$$(\Z[b],\: \Z[b,r]/(r^2+2r))$$
with the right unit
$b\mapsto b+2r$. It can be used to calculate 
the  $E_2$-term of the generalized Adams spectral sequence
$${}^{Y^\prime} E^{s,t}_r(ko)\Rightarrow ko_{t-s}\ .$$
The inclusion $U(2)\to U(4)$ induces a map of ring spectra
$Y^\prime\to Y$, and therefore a morphism of spectral sequences
${}^{Y^\prime} E^{*,*}_*(ko)\to {}^{Y}  E^{*,*}_*(ko)$. 
One can check that this map is an isomorphism from the $E_2$-term on.
The $E_1$-term is the reduced  cobar complex
of the $Y^\prime_*Y^\prime$-comodule $Y^\prime_*ko$.
It consists of torsion-free abelian groups. The same is true after base-change to $Y$.

Using the chain of invariant ideals
$$(2)\subseteq (2,b)\subseteq \Z[b,c]/(r^2+2r)$$
and a routine calculation with Bockstein spectral sequences we obtain  
$${}^YE_2^{*,*}(ko)\cong \Z[b,r]/(2r)\ ,\quad |b|=(0,4)\ ,\: |r|=(1,2)\ .$$
 The equality $ko_3=0$ forces the differential $d_3b=r^3$
which determines $d_3$ algebraically. By sparseness, the spectral sequence degenerates
at $E_4$. Its $E_4=E_\infty$-term reproduces the familiar computation of $ko_*$.
The result is visualized in the following chart. \newpage
 \centerline{\includegraphics[scale=0.8]{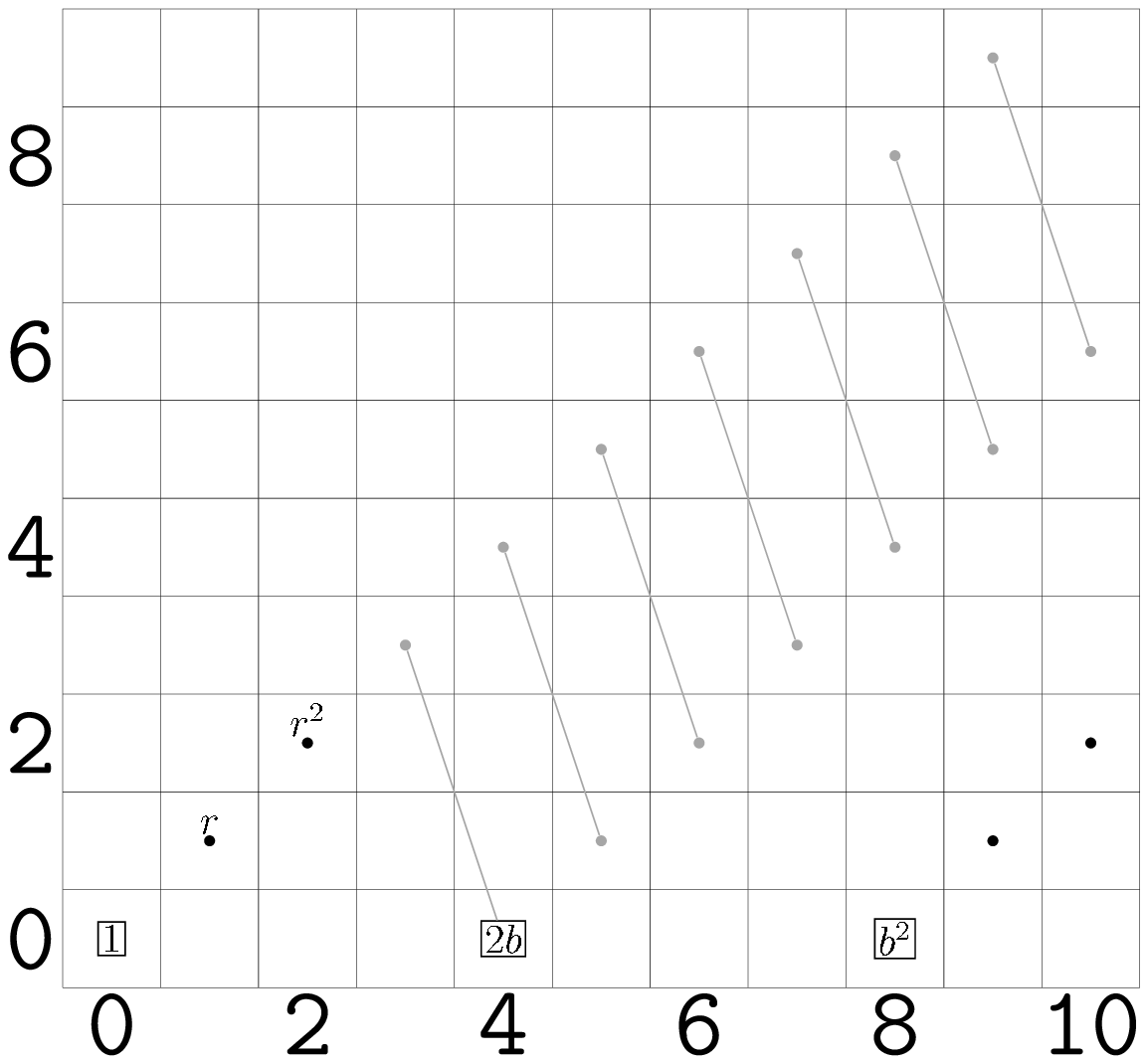}}
\mbox{}\vspace{1cm}\\\centerline{\scriptsize${}^YE_r^{s,t}(ko)\Rightarrow ko_*$}\vspace{1cm}
The horizontal (resp. vertical) axis represents $t-s$ (rep. $s$).
Boxes indicate copies of $\Z$ and dots copies of $\Z/2\Z$. The $E_\infty$-term
is black. The $E_2$-term and the differentials are gray.
 The element $b^2\in {}^YE^{0,8}_2(ko)$ is a periodicity  generator  of dimension $8$.
 
Now we turn to $ko\Q/\Z$.
We have a distinguished triangle
\begin{equation}\label{hdzqdwqdqwdqwd}
ko\to ko\Q\to ko\Q/\Z\stackrel{l}{\to} \Sigma ko
\end{equation}
in $ko$-modules.
Since $Y_*ko$ is torsion-free the map $Y_*(l)$ vanishes. The geometric boundary theorem applies and shows that
the second map in the  exact sequence
$$0\to {}^Y E_2^{s,t}(ko)\otimes \Q/\Z\to {}^Y E^{s,t}_2(ko\Q/\Z)\stackrel{\partial}{\to} {}^Y E_2^{s+1,t}(ko)_{tors}\to 0$$
is compatible with the differential $d_3$.
We get the following chart of ${}^Y E^{*,*}_*(ko\Q/\Z)$.

\vspace{1cm}
\centerline{\includegraphics[scale=0.8]{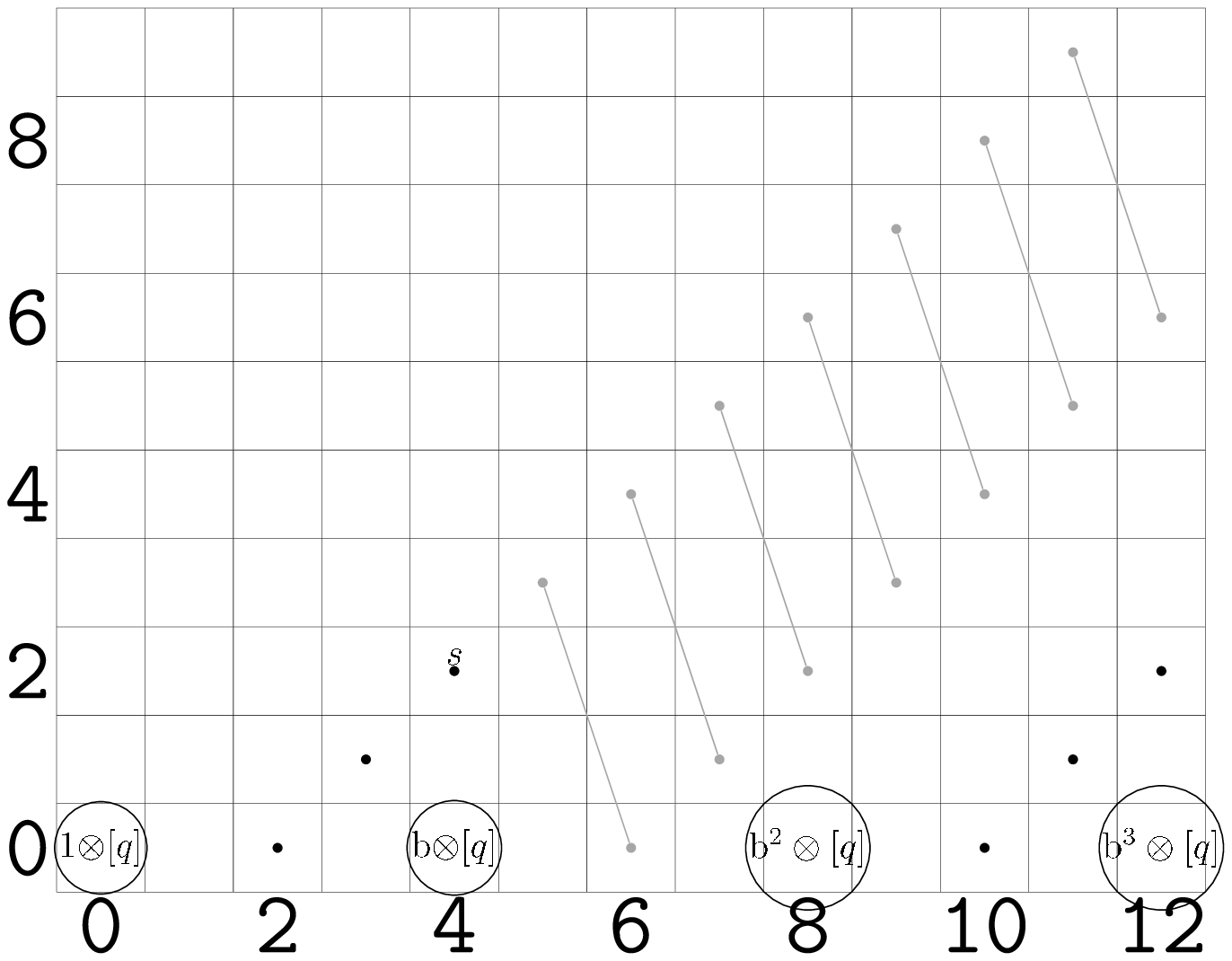}}
\vspace{1cm}
\centerline{\scriptsize${}^YE_r^{s,t}(ko\Q/\Z)\Rightarrow ko\Q/\Z_*$}\vspace{1cm}
The circles indicate copies of $\Q/\Z$ whose elements can be written as $b^k\otimes [q]$ with  $q\in \Q$ and $[q]\in \Q/\Z$. The dots indicate copies of $\Z/2\Z$.
All indicated differentials are forced by the compatibility of $\partial$ with $d_3$. The spectral sequence degenerates at $E_4=E_\infty$. Starting in dimension $2$, the $E_\infty$-term is $8$-periodic with periodicity generator $b^2$ acting via the module structure under ${}^YE^{*,*}_{\infty}(ko)$.

For our applications we are especially interested 
in $ko\Q/\Z_{4m}$.
For all $m\ge 0$ we have
$${}^Y E_2^{0,4m}(ko\Q/\Z)\cong \Q/\Z\ , \quad {}^Y E_2^{2,8m+6}(ko\Q/\Z)\cong \Z/2\Z\ .$$
We let $s\in {}^Y E_2^{2,6}(ko\Q/\Z)\cong \Z/2\Z$ be the non-trivial element. It is a permanent cycle. 
There is an extension problem in dimensions $8m-4$. By the $8$-periodicity it suffices to solve it in dimension $4$.
The long exact sequence in homotopy groups associated with
the triangle (\ref{hdzqdwqdqwdqwd}) gives isomorphisms
$$ko\Q/\Z_{4}\cong ko_{4}\otimes \Q/\Z\cong \Q/\Z\ .$$
This solves the extension problem for $ko\Q/\Z_{4}$ as follows:  The Adams filtration of 
$ko\Q/\Z_{4}\cong  \Q/\Z$ is given by 
$\Z/2\Z\subset \Q/\Z$, where the subgroup is in filtration $\ge 2$.
 We let $\vartheta\in \Filt^2  ko\Q/\Z_{4}$ denote the non-trivial element. It is  detected by 
$[\vartheta]_2=s$. 
We now observe from the chart of the spectral sequence ${}^Y E_2^{*,*} (ko\Q/\Z)$
that $$\Filt^{3} ko\Q/\Z_{*}=0\ .$$
In the following we introduce some  notation to write elements in $ko\Q/\Z_*$ and their
relatives in ${}^YE_\infty^{*,*}(ko\Q/\Z)$.
Let $a_m\in ko_{8m-4}\cong \Z$ be the generator detected by the uniquely determined element $[a_m]_2=2b^{2m-1}$. For $q\in \Q$ we can then consider the class $[q a_m]\in ko\Q/\Z_{8m-4}\cong ko\Q_{8m-4}/ko_{8m-4}$.
If $2q\not\in \Z$, then this class is detected by $2b^{2m-1}\otimes [q]\in
{}^Y E_2^{0,8m-4}(ko\Q/\Z)$. If $q=\frac{1}{2}$, then 
$[q a_m]$ is detected by the unique element $[qa_m]_2=b^{2m-2}s\in {}^Y E_2^{2,8m-2}(ko\Q/\Z)$.

If $a\in ko_{8m}\cong \Z$ is detected by $[a]_2=b^{2m}\in {}^YE^{0,8m}_2(ko)$,
then for $q\in \Q$ the class $[qa]\in ko\Q/\Z_{8m}$ is detected by
$b^{2m}\otimes [q]\in {}^Y E_2^{0,8m}(ko\Q/\Z)$.

It is easy to see that for a spectrum $X\in \{ko,\ ko\Q/\Z\}$ and any prime $p$ we have
$${}^Y E_*^{*,*}(X_{(p)})\cong {}^Y E_*^{*,*}(X)\otimes \Z_{(p)}\ ,$$
where $X_{(p)}$ denotes the $p$-localization of the spectrum $X$.
The above spectral sequences   can be localized at $p=2$ without any essential changes.
If localized at a prime $p\not=2$ they simplify considerably.
In this  case the spectral sequences degenerate at the $E_2$-terms,  and
$${}^Y E^{*,*}_2(ko_{(p)} )\cong \Z_{(p)}[b]\ ,$$
$${}^Y E_2^{0,4m}(ko_{(p)}\Q/\Z)\cong \Q/\Z_{(p)}$$
and vanishes elsewhere.

\section{Calculations at $2$}

\subsection{The result}\label{diwedqwdqwd}

In this section we analyze the   localization at $2$
$$b^{tmf}:tmf_{(2),4m-1}\to T_{(2),2m}\ ,\quad  m\ge 1\ ,$$
of the maps $b^{tmf}:tmf_{4m-1}\to T_{2m}$ constructed in Subsection \ref{uiffwefwefewf},
where $$T_{(2),2m}:=(T_{2m})_{(2)}:=T_{2m}\otimes \Z_{(2)}$$ is the localization of $T_{2m}$ at $2$.
 In view of the surjectivity of the Witten genus
$\sigma:MString_{*}\to tmf_{*}$ (\cite[Theorem 6.25]{MR1989190}) and Proposition \ref{udiqwdqwdqwdq} the result implies 
a non-triviality statement about the homomorphism $b^{top}:A_{4m-1}\to T_{2m}$
for an infinite number of  $m\ge 1$.

In order to state the result we must name some elements of $tmf_{(2),*}$.
To this end we  first recall some parts of the calculations of \cite{MR2508200} of the spectral sequence ${}^Y E_*^{*,*}(tmf_{(2)})$ and of $tmf_{(2),*}$.
Let $\Delta\in \cM^\Z_{12}$ be the unique cusp form which is normalized such that $\Delta=q+\dots$.
There exists an unique element 
$\Delta_2\in {}^Y E_2^{0,24}(tmf_{(2)})\cong \Z^2_{(2)}$
which corresponds to $\Delta$ under the identification (\ref{ghagdastdtztzqw}).
We know that  $8\Delta_2$ is a permanent cycle.
 The element $\Delta_2^8\in   {}^Y E_2^{0,192}(tmf_{(2)})$ is a permanent cycle. It detects a unique element 
 $\Delta^8\in tmf_{(2),192}$ which is a periodicity  generator of $tmf_{(2),*}$.

The map $\sigma : MString\to tmf$ induces an isomorphism $\sigma : MString_3\to tmf_3$.
The generator $g\in  MString_3$ defined in (\ref{reuwrewrwer}) maps to a generator
$\nu:=\sigma(g)$ of $tmf_{(2),3}$, where we omit   the localization map
$tmf\to tmf_{(2)}$ from notation. The  order of $\nu$ is $8$,  it belongs to  filtration $1$, and it is detected by  a unique element $[\nu]_2\in {}^Y E_2^{1,4}(tmf_{(2)})$.
Let $\eta\in tmf_{1}$ be the image of the Hopf map in $S_1$ under the unit $S\to tmf$.
It is detected by $[\eta]_2\in {}^Y E_2^{1,2}(tmf_{(2)})$.
For $i\in \{14,38,74,110,134\}$  we let  $[a_i]_2\in {}^Y E_2^{1,i}(tmf_{(2)})\cong  \Z/2\Z$ denote  generators of the corresponding groups.

The upper half of the fourth column of  the following table displays all additive generators of ${}^Y E_\infty^{1,4m}(tmf_{(2)})\subseteq E_2^{1,4m}(tmf_{(2)})$ and  dimension $4m-1<195$.  We fix elements of
$tmf_{(2),4m-1}$ detected by these generators and which are named in the second  column.
These elements are unique up to elements of filtration $\ge 3$ (with the exception of $\nu$ which we have fixed above). 
The third  column lists the order of the elements in the first column. If we multiply the elements of order $8$ (or $4$, respectively) by $4$ (or $2$, respectively), then we get elements in filtration $\ge 3$.
These are detected by elements in ${}^YE_\infty^{3,4m}(tmf_{(2)})$ which do not appear separately in the table.
The lower half  of the fourth column lists the remaining permanent cycles of ${}^Y E_2^{3,4m+3}(tmf_{(2)})$ in dimension $4m-1<192$.

The second column is a complete list of additive generators of
$$\Filt^{ 1}tmf_{(2),4m-1}/\Filt^{ 4}tmf_{(2),4m-1}\ .$$
$$\begin{array}{|c|c|c|c|c|c|}\hline 
m&name &  \ord & {}^YE^{*,*}_2(tmf_{(2)})   &b^{tmf}(\dots)      & c\in \{\dots\}  \\ \hline\hline
1 &\nu  &8                    &[\nu]_2              &[\frac{3}{8}]&  \\  \hline
7&2\nu \Delta&4 &2[\nu]_2\Delta_2&[\frac{c}{4}\Delta]  &1,3\\\hline
13&\nu \Delta^2&8 &[\nu]_2\Delta^2_2&[\frac{c}{8}\Delta^2] & 1,5 \\\hline
25&\nu \Delta^4&8 &[\nu]_2\Delta^4_2&[\frac{c}{8}\Delta^4] & 1,5 \\\hline
31&2\nu \Delta^5&4 &2[\nu]_2\Delta^5_2&[\frac{c}{4}\Delta^5] & 1,3 \\\hline
37&\nu \Delta^6&8 &[\nu]_2\Delta^6_2&[\frac{c}{8}\Delta^6] &1,5 \\\hline\hline
 4&\eta a_{14}&2&[\eta]_2 [a_{14}]_2&0&\\\hline10&\eta a_{38}&2&[\eta]_2 [a_{38}]_2&0&\\\hline19&\eta a_{74}&2&[\eta]_2 [a_{74}]_2&0&\\\hline29&\eta a_{110}&2&[\eta]_2 [a_{110}]_2&0&\\\hline34&\eta a_{134}&2&[\eta]_2 [a_{134}]_2&0&\\\hline
 \end{array}\ .
$$
The table can be continued in the obvious $192$-periodic way using the multiplication by $\Delta^8$.
There are more elements in dimensions $4m-1$ which are products of
the listed elements with other elements in dimensions divisible by $4$. All these products 
are of filtration $\ge 4$. 
The main result of the present section is the following proposition.
\begin{prop}\label{uidodqwdqwd}
\begin{enumerate}
\item The values of $b^{tmf}$ on the elements listed in the second column are represented by the formal power series listed in the fifth  column of the table. 
\item The map $b^{tmf}$ annihilates all elements in filtration $\ge 4$.
\item\label{udiqwdhqwdqwdwqd}  For all $m\ge 0$ and $[m]_{\mbox{mod $48$}}\in \{1,7,13,25,31,37\}$
the map $b^{tmf}$ induces an injective map
$$\bar b^{tmf}:\Filt^{ 1}tmf_{(2),4m-1}/\Filt^{4}tmf_{(2),4m-1}\to T_{(2),2m}\ .$$
\item For all $m\ge 0$ not listed in \ref{udiqwdhqwdqwdwqd}
the map $\bar b^{tmf}:tmf_{(2),4m-1}\to T_{(2),2m}$ is trivial. 
 \end{enumerate}
\end{prop}
The remainder of the present section is devoted to the proof of this proposition.

\subsection{Transition to $\Q/\Z$-theory}\label{fuweifwefweftztwf}

 From Subsection  \ref{udiqwdwqd} 
we have an exact sequence
\begin{equation}\label{u7diqwdqwd444}
0\to {}^Y E_2^{s,t}(tmf_{(2)})\otimes \Q/\Z\to {}^Y E_2^{s,t}(tmf_{(2)}\Q/\Z)\stackrel{\partial}{\to} {}^Y E_2^{s+1,t}(tmf_{(2)})_{tors}\to 0
\end{equation}
for all $t,s\in \Z$.
The group 
${}^Y E_2^{s,t}(tmf_{(2)})$    is  finite for all  $s\ge 1$ and $t\in \Z$.
This implies that
$ {}^Y E_2^{s,t}(tmf_{(2)})\otimes \Q/\Z=0$
for $s\ge 1$.
Note that 
${}^Y E_2^{0,t}(tmf_{(2)})\otimes \Q/\Z$ is a divisible group.
There are no non-trivial maps from divisible groups to finite  groups. Hence
the  elements in the subgroup   $ {}^Y E_2^{0,t}(tmf_{(2)})\otimes \Q/\Z\subseteq {}^Y E_2^{0,t}(tmf_{(2)}\Q/\Z)$
are  annihilated by all differentials, i.e. they are permanent cycles.
In particular we get an embedding
$$0\to {}^Y E_2^{0,t}(tmf_{(2)})\otimes \Q/\Z\to  {}^Y E_\infty^{0,t}(tmf_{(2)}\Q/\Z)\ .$$
By an inspection of the calculation  \cite{MR2508200} we see 
that the   maps
$$\delta_r: {}^Y E_r^{s,4m}(tmf_{(2)}\Q/\Z)\to  {}^YE_r^{s+1,4m}(tmf_{(2)})$$
 are isomorphisms for all $r\ge 2$ and $s\ge 1$.
Indeed, they map cycles to cycles.
Since the differentials of ${}^Y E_r^{*,*}(tmf_{(2)}\Q/\Z)$ annihilate the subgroups ${}^Y E_2^{0,t}(tmf_{(2)})\otimes \Q/\Z$ and otherwise are induced by the differentials of
${}^Y E_r^{*,*}(tmf_{(2)})$
 we further conclude that if  $\delta_r(x)$ is a boundary, then so is
$x$.
We therefore  have isomorphisms
$$ \delta_r: \frac{ {}^Y E_r^{0,4m}(tmf_{(2)}\Q/\Z)}{{}^Y E_2^{0,4m}(tmf_{(2)})\otimes \Q/\Z}\stackrel{\sim}{\to}  E_r^{1,4m}(tmf_{(2)})$$
for all $r\ge 2$ and $r=\infty$.

\subsection{The map to $ko[[q]]$}\label{uwiefwefwefwefewfwfef}

The morphism of ring spectra $w:tmf\to ko[[q]]$ induces a morphism $E_r(w)$ of $Y$-based  generalized Adams spectral sequences.
Note that
${}^Y E_2^{*,*}(ko[[q]])\cong {}^Y E^{*,*}_2(ko)[[q]]$.
The spectra $ko[[q]]_{(2)}$, $tmf_{(2)}\Q/\Z$ and $ko[[q]]_{(2)}\Q/\Z$   are modules over $tmf_{(2)}$.
We get induced module structures on the spectral sequences.

Recall the element $\Delta_2\in {}^Y E_2^{0,24}(tmf_{(2)})$ introduced in Subsection  \ref{diwedqwdqwd}.
We have 
$$E_2(w)(\Delta_2)=b^{6} \Delta\in {}^Y E_2^{0,24}(ko[[q]]_{(2)})\cong {}^Y E_2^{0,24}(ko)[[q]]_{(2)}\ .$$

We let $\bar w=w\wedge \id_{M\Q/\Z}:tmf_{(2)}\Q/\Z\to ko[[q]]_{(2)}\Q/\Z$ denote the map induced by $w$.
Since $$\Filt^{ 3}ko[[q]]_{(2)}\Q/\Z_*=0$$ it is clear that
$\bar w$ annihilates $\Filt^{3}tmf_{(2)}\Q/\Z_*$. It now follows from Lemma \ref{duqwiduqwdqwdqwd} that $b^{tmf}$ annihilates
$\Filt^{4}(tmf_{(2)})$.  This is Assertion 2. of Proposition \ref{uidodqwdqwd}.
Consequently $b^{tmf}$ annihilates all elements except for those listed in the table in Subsection \ref{diwedqwdqwd} in dimensions $4m-1$.

  We now explain how we get the listed  values of $b^{tmf}$.

Recall the generator $\nu\in tmf_{(2),3}$ of order $8$. We choose a  lift
$\nu_{\Q/\Z}\in tmf_{(2)}\Q/\Z_{4}$ as in (\ref{diagr2}).
By (\ref{dzuwdzquwd}) know that
$b^{top}(g)=[\frac{1}{24}]\in T_2= \R[[q]]/\Z[[q]]$. Therefore, using that
$[\frac{1}{24}]=[\frac{3}{8}]$ in $\Q/\Z_{(2)}$, we get by Lemma \ref{duqwiduqwdqwdqwd}
$$\bar w( \nu_{\Q/\Z})=[-\frac{3}{8}] \in  ko[[q]]_{(2)}\Q/\Z_4\cong \Z_{(2)}[
[q]]\otimes \Q/\Z\ .$$
From the discussion at the end of Subsection \ref{uiqdqwdqwdqwd234324234} it follows that
$$E_2(\bar w)([ \nu_{\Q/\Z}]_2)
=b\otimes [-\frac{3}{4}]\in E_2^{0,4}(ko[[q]]_{(2)}\Q/\Z)\cong \Z[[q]]_{(2)}\otimes \Q/\Z\ .$$

Let $\lambda=2\nu \Delta\in \Filt^{1}tmf_{(2),27}$ be detected by $[\lambda]_2:=2[\nu]_2\Delta_2$.
Then we can choose a lift
$  \lambda_{\Q/\Z}\in tmf_{(2)}\Q/\Z_{27}$ which is detected by 
$[\lambda_{\Q/\Z}]_2=2[ \nu_{\Q/\Z}]_2 \Delta_2$ (note that the lift $  \lambda_{\Q/\Z} $ is unique up to elements coming from $tmf_{(2)}\Q_{28}\cong \Q$).
This follows from the condition  $\partial [ \lambda_{\Q/\Z}]_2=[\lambda]_2$\ ,
where $\partial$ is as in (\ref{u7diqwdqwd444}), and
$\partial (2[\nu_{\Q/\Z}]_2 \Delta_2)=2\Delta_2 \partial ([ \nu_{\Q/\Z}]_2)$.

  We conclude that
$$E_2(\bar w)( \lambda_{\Q/\Z})=E_2(\bar w)([ \nu_{\Q/\Z}]_2) 
2\Delta_2=b^7\otimes [\frac{1}{2}\Delta]\in E_2^{0,28}(ko[[q]]_{(2)}\Q/\Z)\ .$$
Lemma \ref{duqwiduqwdqwdqwd} implies
that
$$b^{tmf}(2\nu \Delta)=[\frac{c}{4}\Delta]\in T_{(2),28}\ ,$$
where $c=1$ or $3$. This indeterminacy is unavoidable since 
$2\nu \Delta$ is only determined up to elements of higher filtration.

By a similar argument 
$$b^{tmf}(\nu\Delta^2)=[\frac{c}{8}\Delta^2]\in T_{(2),52}$$
with $c\in \{3,7\}$. In a similar manner we get the value of  $b^{tmf}$ on the remaining entries of the upper part of the second column. 

Turning to the elements listed in the lower part of the second column one 
shows $b^{tmf}(\eta a_i)=0$ as follows.
Note that $\eta a_i\in Filt^3 tmf_{(2),i+1}$.
Therefore $b^{tmf}(\eta a_i)$ would be represented by an element in
$Filt^{2} ko[[q]]\Q/\Z_{i+2}$.
In all cases we have $i+2\equiv 0 \:\mbox{mod $8$}$. 
But then 
$Filt^2 ko[[q]]\Q/\Z_{i+2}=0$.
We thus have shown Assertions 1 and 4 of Proposition \ref{uidodqwdqwd}.

It remains to show Assertion 3. It suffices to show that
$\ord(b^{tmf}(x))=\ord(x)$ for the generators of the cyclic groups
$\Filt^{1}tmf_{(2),4m-1}/\Filt^4tmf_{(2),4m-1}$ listed in the first column
of the table and their multiples by $\Delta^8$. The assertion immediately follows from Corollary \ref{zdwdqwdqwdwddqwdq}. 
\hB

\section{Calculations at $3$}

\subsection{The result}\label{diwedqwdqwd1}

Compared with the localization at $2$, the structure of the spectral sequence ${}^Y E^{*,*}_*(tmf_{(3)})$ and of $tmf_{(3),*}$
is much simpler. We again start by recalling 
the calculations of \cite{MR2508200}.
There is an unique element $\Delta_2\in {}^Y E_2^{0,24}(tmf_{(3)})$
which corresponds to $\Delta$ under the identification (\ref{ghagdastdtztzqw}).
We know that 
 $3\Delta_2$ is a permanent cycle. 
The element $\Delta_2^3\in  {}^Y E_2^{0,72}(tmf_{(3)})$ is a permanent cycle and detects a unique element
 denoted by $\Delta^3\in tmf_{(3),72}$ which is  a periodicity generator. 
 
The generator $g\in  MString_3$ defined in (\ref{reuwrewrwer}) maps to a generator
$\nu:=\sigma(g)$ of $tmf_{(3),3}$ of order $3$ in filtration $1$ detected by  $[\nu]_2\in {}^Y E_2^{1,4}(tmf_{(3)})$.
 The following table gives the complete list of additive generators of $tmf_{(3),4m-1}$ for 
 $4k-1< 75 $.  All these elements live in filtration one  and   are  uniquely determined.
The table can be continued in the obvious $72$-periodic way using multiplication by $\Delta^3$.
$$\begin{array}{|c|c|c|c|c|c|}\hline 
m&name &  \ord  &{}^YE^{*,*}_2(tmf_{(3)})   &b^{tmf}(\dots)        \\ \hline\hline
1 &\nu  &3                  &[\nu]_2              &[\frac{2}{3}]  \\  \hline
7&\nu \Delta&3 &[\nu]_2\Delta_2&[\frac{2}{3}\Delta]    
 \\\hline
  \end{array}\ .
$$
The third column displays the order of the element.  The last two columns display the elements in ${}^YE^{*,*}_2(tmf_{(3)})$ which detect the elements in the second column, and the value
of the invariant $b^{tmf}$ in $T_{(3),2k}$,  expressed  as a class of a rational formal power series.

\begin{prop}\label{zwdqdwqd}
\begin{enumerate}
\item The value of $b^{tmf}$ on the generator listed in the second column of the table is as indicated in the last column of the table.
\item For all $m\ge 1$ we have injections
$$b^{tmf}:tmf_{(3),4m-1}\to T_{(3),2m}\ .$$
\end{enumerate}
\end{prop}

The following two subsections are devoted to the proof of this proposition.

\subsection{Transition to $\Q/\Z$-theory}

 By Subsection \ref{udiqwdwqd}
we have an exact sequence
\begin{equation}\label{u7diqwdqwd4443}
0\to {}^Y E_2^{s,t}(tmf_{(3)})\otimes \Q/\Z\to {}^Y E_2^{s,t}(tmf_{(3)}\Q/\Z)\stackrel{\partial}{\to} {}^Y E_2^{s+1,t}(tmf_{(3)})_{tors}\to 0\ .
\end{equation}
As in Subsection \ref{fuweifwefweftztwf}  we see that 
the  elements in the subgroup   $ {}^Y E_2^{0,t}(tmf_{(3)})\otimes \Q/\Z\subseteq {}^Y E_2^{0,t}(tmf_{(3)}\Q/\Z)$
are  permanent cycles and
  we get an embedding
$$0\to {}^Y E_2^{0,t}(tmf_{(3)})\otimes \Q/\Z\to  {}^Y E_\infty^{0,t}(tmf_{(3)}\Q/\Z)\ .$$
By an inspection of the calculation  \cite{MR2508200} we see  again
that the  maps
$$\delta_r: {}^Y E_r^{s,4m}(tmf_{(3)}\Q/\Z)\to  E_r^{s+1,4m}(tmf_{(3)})$$
 are isomorphisms for all $r\ge 2$ and $s\ge 1$, so that the maps 
$$ \delta_r: \frac{ {}^Y E_r^{0,4m}(tmf_{(3)}\Q/\Z)}{{}^Y E_2^{0,4m}(tmf_{(3)})\otimes \Q/\Z}\stackrel{\sim}{\to}  E_r^{1,4m}(tmf_{(3)})\ $$
are isomorphisms for all $r\ge 2$ and $r=\infty$.

\subsection{The map to $ko[[q]]$}\label{qwddwqdqwdoiuqwdoiwqd}

The element  $\Delta_2\in {}^Y E_2^{0,24}(tmf_{(3)})$ satisfies
$$E_2(w)(\Delta_2)=b^{6} \Delta\in {}^Y E_2^{0,24}(ko[[q]]_{(3)})\cong {}^Y E_2^{0,24}(ko)[[q]]_{(3)}\ .$$
Recall the generator $\nu\in tmf_{(3),3}$ of order $3$ and the lift
$ \nu_{\Q/\Z}\in tmf_{(3)}\Q/\Z_{4}$ as in (\ref{diagr2}). 

 By (\ref{dzuwdzquwd}) know that
$b^{top}(g)=[\frac{1}{24}]\in T_2\cong \Q/\Z[[q]]$. Therefore, using that
$[\frac{1}{24}]=[-\frac{1}{3}]$ in $\Q/\Z_{(3)}$, we get
$$\bar w( \nu_{\Q/\Z})=[\frac{1}{3}] \in  ko[[q]]_{(3)}\Q/\Z_4\cong \Z[
[q]]_{(3)}\otimes \Q/\Z\ .$$
With Lemma  \ref{duqwiduqwdqwdqwd} this gives $$b^{tmf}(\nu)=\left[\frac{2}{3}\right]\in T_{(3),2}\ .$$

Furthermore, we have
$$E_2(\bar w)([\tilde \nu]_2)
=b\otimes [\frac{2}{3}]\in E_2^{0,4}(ko[[q]]_{(3)}\Q/\Z)\cong \Q/\Z[[q]]_{(3)}\ .$$
 
Let $\lambda=\nu \Delta\in \Filt^{1}tmf_{(3),27}$ be detected by $[\lambda]_2:=[\nu]_2\Delta_2$.
Then we can choose a lift
$ \lambda_{\Q/\Z}\in tmf_{(3)}\Q/\Z_{27}$ which is detected by 
$[ \lambda_{\Q/\Z}]_2=[\tilde \nu]_2 \Delta_2$ (note that the lift $ \lambda_{\Q/\Z}$ is unique up to elements coming from $tmf_{(3)}\Q_{28}\cong \Q$).
 
 We conclude that
$$E_2(\bar w)( \lambda_{\Q/\Z})=E_2(\bar w)( \nu_{\Q/\Z}]_2) 
\Delta_2=b^7\otimes [\frac{2}{3}\Delta]\in E_2^{0,28}(ko[[q]]_{(3)}\Q/\Z)\ .$$
This gives $$b^{tmf}(\nu \Delta)=\left[\frac{2}{3}\Delta\right]\in T_{(3),14}\ .$$
This finishes the proof of Assertion 1 of Proposition \ref{zwdqdwqd}. .
 
For Assertion 2.
 it suffices to show that
$\ord(b^{tmf}(x))=\ord(x)$ for the generators of the cyclic groups
$tmf_{(3),4m-1}$ listed in the first column
of the table and their multiples by $\Delta^3$. The assertion immediately follows from Corollary \ref{zdwdqwdqwdwddqwdq}. 
 \hB

\subsection{How to detect $\nu\Delta$?}\label{ziddqwdqwdwq}

The equality $b^{top}=b^{geom}$ can be used to do some explicit calculations of the tmf-valued Witten genus. In the following we discuss an  example.
Let  $(M,\alpha^{top})$ be a $27$-dimensional closed string manifold. It represents a string bordism class $[M,\alpha^{top}]\in MString_{27}$. Moreover,
$\sigma([M,\alpha^{top}])= c \Delta \nu\in  tmf_{(3),27}$ for a uniquely determined 
 $c\in \Z/3\Z$.  
The problem consists in the explicit  calculation of the value $c$.

First observe that  by the calculations in Subsection \ref{diwedqwdqwd1} we have an equality of classes
in $T_{(3),27}$
$$[b^{top}(M,\alpha^{top})]=[\frac{2c}{3} \Delta]\ .$$
Expand $b^{top}(M,\alpha^{top})=b_0+qb_1+\dots$.
Then we form
$$b^{top}(M,\alpha^{top})-b_0 E_{14}=(b_1-24 b_0)q+\dots \ . $$
Note that $3 (b^{top}(M,\alpha^{top})-b_0 E_{14})\in \Z[[q]]_{(3)}$.
We have
$$\frac{2}{3}\Delta= \frac{2}{3} q+\dots \quad \ ,\quad\ .$$
Therefore we get
$$c=-3(b_1-24 b_0)  \: (\modu 3)\ .$$
Let $Z$ be a spin zero bordism of $M$.
Then
specializing Lemma (\ref{zdqwud}) we get
\begin{eqnarray*}
c&=&-\langle\tilde p_1\cup  \left({\frac {1967}{729}}\,{N_{{4}}}^{3}+{\frac {356}{243}}\,N_{{8}}{p_{{1}
}}^{2}+{\frac {2575}{2187}}\,N_{{6}}{p_{{1}}}^{3}\right.\\&&\left.+{\frac {152}{81}}\,
N_{{8}}N_{{4}}+{\frac {941}{729}}\,N_{{4}}p_{{1}}N_{{6}}+{\frac {6232
}{2187}}\,{p_{{1}}}^{6}\right.\\&&\left.+{\frac {898}{729}}\,N_{{4}}{p_{{1}}}^{4}+{
\frac {541}{243}}\,N_{{10}}p_{{1}}+{\frac {623}{729}}\,{N_{{4}}}^{2}{
p_{{1}}}^{2}\right.\\&&\left.+{\frac {457}{729}}\,N_{{12}}+{\frac {2398}{2187}}\,{N_{{
6}}}^{2}\right), [Z,M]\rangle \mod 3
 \end{eqnarray*}
Here $N_{2i}:=N_{2i}(p_1,p_2,...)$ denotes the $i$'th Newton polynomial and we use the abbreviations $p_i:=p_i(TZ)$ and $\tilde p_1:=\tilde p_1(TZ,\alpha^{top})$ for the class   defined in (\ref{uidquwdqwdqwdqwdqwdqwdq}).
It is a non-trivial integrality statement that the right-hand side is an integer (before taking its value modulo $3$).

\section{Geometric and homotopy theoretic string bordism}\label{udidqwdqwdqwdqwdqwdqwd}

In the present paper we use a geometric picture of string bordism.
The goal of  the present subsection is to give the geometric definition of  string bordism and to relate it with the 
homotopy theoretic picture. 
 This material is familiar to geometers and
we include it mainly as a reference.

We first recall some facts about string structures. The following diagram   gives the unstable variation
 of the Postnikov tower (\ref{uifwefwefwefewf}). For $n\ge 3$ we consider
\begin{equation}\label{uqiwdqwdqwdwd}
\xymatrix{&BString(n)\ar[d]&\\ X\ar[ddr]_{V}\ar@{.>}[ur]^{\alpha^{top}}\ar[r]_/-1em/{\xi}&BSpin(n)\ar[d]\ar[r]^{\frac{p_1}{2}}&K(\Z,4)\\
&BSO(n)\ar[r]^{w_2}\ar[d]&K(\Z/2\Z,2)\\&BO(n)\ar[r]^{w_1}&K(\Z/2\Z,1)}\ .
\end{equation}
The space $BSO(n)$ is the homotopy fiber of the first Stiefel-Whitney class $w_1$. Similarly, $BSpin(n)$ and $BString(n)$ are the homotopy fibers of the second Stiefel-Whitney class $w_2$ and the spin refinement $\frac{p_1}{2}$ of the first Pontrjagin class, respectively.

Let  $V\to X$ be a $n$-dimensional real vector bundle over a  $CW$-complex $X$. We use same symbol in order to denote the classifying map  $V:X\to BO(n)$. An orientation and  spin structure on $V$ is given by the choice of a lift $\xi:X\to BSpin(n)$. A spin bundle is by definition a pair $(V,\xi)$. We will usually drop the $\xi$ from the notation and refer to  $V$ as a spin bundle.
\begin{ddd}
A topological string structure on a spin bundle $V$ is defined as a  homotopy class of lifts  $\alpha^{top}$. The pair $(V,\alpha^{top})$ is called a string bundle.
\end{ddd}
The set of topological string structures on the spin bundle $V$ is a torsor under $H^3(X;\Z)$. We will write the action of $x\in H^3(X;\Z)$ as 
\begin{equation}\label{ufiwefwefwefw}
(x,\alpha^{top})\mapsto \alpha^{top}+x\ . 
\end{equation}

In differential geometry one works with orientations, spin, or string structures on the (unstabilized)
tangent bundle, while in homotopy theory one considers those structures on the stable normal bundle.
In order to transfer structures back and forth between the tangent and stable normal bundles
the two-out-of three principle Lemma \ref{ioqwdqwdqwdwqd} below is crucial.

In the following we describe an abstract setting for this  principle which can inductively  be applied to the layers of the tower  (\ref{uqiwdqwdqwdwd}).
Assume that we have a sequence of maps
$(B(n)\to BO(n))_{n\ge n_0}$ together with homotopy commutative diagrams
$$\xymatrix{B(n)\times B(m)\ar[d]\ar[r]^{\kappa_B}&B(n+m)\ar[d]\\BO(n)\times BO(m)\ar[r]^/0.5em/{\oplus}&BO(n+m)}$$
satisfying the obvious associativity constraints. We further assume that the spaces $B(n)$ are $l-1$-connected for some integer $l$,  and that we have 
a collection of maps $(c_n:B(n)\to K(\pi,l))_{n\ge n_0}$
inducing  isomorphisms in $\pi_l(B(n))\stackrel{\sim}{\to}\pi$ for all $n\ge n_0$ and some fixed abelian group $\pi$.
We define the  $l$-connected spaces $C(n)$ as homotopy  pull-backs
$$\xymatrix{C(n)\ar[r]\ar[d]&\mbox{}*\ar[d]\\
B(n)\ar[r]^{c_n}&K(l,\pi)}\ .$$
Finally we assume that the collection of maps $(c_n)_{n\ge n_0}$ is additive in the sense that the lower squares of the following diagrams are homotopy commutative:
$$\xymatrix{C(n)\times C(m)\ar@{.>}[drr]^{\kappa_C}\ar[ddd]&&&\\&&C(n+m)\ar[ddd]&\\&K(\pi,l)\times K(\pi,l)\ar[rrd]^+&&\\B(n)\times B(m)\ar[rrd]^{\kappa_B}\ar[ru]^{c_n\times c_m}&&&K(\pi,l)\\&&B(n+m)\ar[ur]_{c_{n+m}}&}\ .
$$
Then we get  unique (up to homotopy) lifts $\kappa_C$ of $\kappa_B$ indicated by the dotted arrow which again satisfy
an associativity constraint.

Let $X$ be a $CW$-complex. A real $n$-dimensional vector bundle $V\to X$ with a $B$-structure is a diagram
$$\xymatrix{&B(n)\ar[d]\\X\ar[ur]^\xi\ar[r]^V&BO(n)}\ ,$$
where $V$ is the classifying map of the vector bundle, and the lift $\xi$ defines the $B(n)$-structure.
The pair $(V,\xi)$ will be called a $B$-bundle.
A $C$-structure on a $B$-bundle $(V,\xi)$ is given by a  further lift $\alpha$ as indicated in 
$$\xymatrix{&C(n)\ar[d]\\&B(n)\ar[d]\\X\ar@{.>}[uur]^\alpha\ar[ur]^\xi\ar[r]^V&BO(n)}\ .$$
The set of isomorphism classes of $C$-structures on the $B$-bundle $(V,\xi)$ 
is in one-to-one correspondence with the set of homotopy classes of lifts $\eta$ of $\xi$.
By obstruction theory this set is a torsor under
$H^l(X;\pi)$. If $x\in H^l(X;\pi)$, then we let $\alpha+x$ denote the result of the action of $x$ on $\alpha$.
The maps $\kappa_B $  allow to define the sum of  $B$-bundles 
$$(V,\xi)\oplus (W,\zeta):=(V\oplus W,\kappa_B(\xi,\zeta))\ .$$
Furthermore, using the maps $\kappa_C$ we define
the sum of $B$-bundles with $C$-structures by 
$$((V,\xi),\alpha)\oplus ((W,\zeta),\beta):=((V\oplus W,\kappa_B(\xi,\zeta)),\kappa_C(\alpha,\beta))\ .$$
 One can check that
$\kappa_C$ is bilinear with respect to the action of $H^l(X;\pi)$, i.e.
we have
$$\kappa_C(\alpha+x,\beta)=\kappa_C(\alpha,\beta+x)=\kappa_C(\alpha,\beta)+x\ , \quad \forall x\in H^l(X;\pi)\ .$$
We can now formulate the two-out-of-three principle for $C$-structures on $B$-bundles.
We consider $B$-bundles $(V,\xi)$ and $(W,\zeta)$. 
\begin{lem}
The choice of two $C$-structures out of $\alpha,\beta,\gamma$ on  the $B$-bundles  $$(V,\xi)\ ,\quad (W,\zeta)\ ,\quad (V\oplus W,\kappa_B(\xi,\zeta))$$ uniquely fixes
the third such that  there is an isomorphism of $C$-bundles
$$((V,\xi),\alpha)\oplus((W,\zeta),\beta)\cong ((V\oplus W,\kappa_B(\xi,\zeta)),\gamma)\ .$$
\end{lem}
The characteristic classes which define the layers of the Postnikov tower (\ref{uqiwdqwdqwdwd})
are additive. Hence 
the theory described above can be applied.
 Recalling the stability statements $\pi_0(O(n))=\Z/2\Z$ for $n\ge 1$,
$\pi_1(SO(n))=\Z/2\Z$ and $\pi_3(Spin(n))=\Z$ for 
$n\ge 3$ we thus obtain the following.
\begin{kor}\label{ioqwdqwdqwdwqd}
We have a two-out-of-three principle for 
\begin{itemize}
\item orientations on real vector bundles
\item spin-structures on oriented vector bundles of rank $\ge 3$ and
\item string structures on spin bundles of rank $\ge 3$.
\end{itemize} 
\end{kor}
In contrast,  this principle does not apply to (homotopy classes of) framings   as is obvious from
$TS^2\oplus\R\simeq S^2\times\R^3$.

We now describe the geometric picture of string bordism.
A spin manifold $M$ is a manifold with a chosen orientation and spin structure on its tangent bundle $TM$.
\begin{ddd}\label{uidqwdqwdqwdqwd}
A string manifold $(M,\alpha^{top})$ is a spin manifold $M$ with a chosen topological string structure $\alpha^{top}$ on the spin bundle  $TM$.
\end{ddd}

If $(Z,\beta^{top})$ is a string manifold with boundary $M$, then $M$
has an induced string structure $\alpha^{top}$ defined as follows.
The inner normal field induces a decomposition
\begin{equation}\label{udiqduqwdwqd}TZ_{|M}\cong TM\oplus \Theta_\R\ ,
\end{equation} where $\Theta_\R:= M\times \R$
denotes the trivial bundle. A trivial bundle has a canonical orientation, as well as a canonical spin structure, and a canonical string structure $can$. We define the orientation, spin structure, and  string structure on $TM$ by the two-out-of-three principle such that (\ref{udiqduqwdwqd})
becomes a decomposition of structured bundles.

We can now define the $n$'th geometric string bordism homology group $\Omega^{String}_n(X)$ of a space $X$ in the usual manner by cycles and relations. A cycle for $\Omega^{String}_n(X)$ is a pair $((M,\alpha^{top}),f)$ of an $n$-dimensional closed string manifold $(M,\alpha^{top})$ and a map $f:M\to X$. The set of isomorphism classes of cycles forms a semigroup with respect to disjoint union.  The group $\Omega^{String}_n(X)$ is the defined as the quotient of the group completion of this semigroup by the group generated by string manifolds with maps to $X$ which are boundaries of $n+1$-dimensional compact string manifolds with maps to $X$. The remaining structures for the homology theory
(e.g. the suspension isomorphism) are defined in the standard way and will be neglected in the discussion below, too.

We have the following consequence of the Pontrjagin-Thom construction and the
two-out-of-three principle Lemma \ref{ioqwdqwdqwdwqd}. 
\begin{lem}\label{udiqdqwdwqd}
There is a canonical isomorphism 
$$\Omega^{String}_*(X)\cong MString_*(X)\ $$
for all $*\ge 3$.
\end{lem}
\proof
The standard Thom-Pontrjagin construction produces an isomorphism between the homotopy theoretic string bordism and a version of the geometric string bordism
with topological string structures on stable normal bundles \cite[Theorem, Chapter II]{stongcob}. Using the two-out-of-three principle Lemma \ref{ioqwdqwdqwdwqd} we can go back and forth between topological string structures on stable normal bundles and on the (unstable) tangent bundle. \hB

 In the main body of the text, we will use the notation $MString_*(X)$ in order to denote both, the geometric and the homotopy theoretic string bordism groups of $X$.

\section{Modular forms and the structure of $T_{2m}$}\label{fzuefwefwefewf}
\newcommand{\MMM}{{\mathfrak{M}}}
\newcommand{\WWW}{{\mathfrak{W}}}
\newcommand{\CCC}{{\mathfrak{C}}}
The present subsection has two goals. One is to make the structure of the group $T_{2m}$ more explicit. In applications the  important problem is to decide whether a formal power series
$f\in \R[[q]]$ represents a non-trivial element in $T_{2m}$.  This can be solved using Lemma \ref{invariant}.
The other goal is to fix the isomorphism between the $E_2$-term ${}^YE_2^{0,4m}(tmf)$ of the $Y$-generalized based Adams sequence for $tmf$ and the space of  modular forms $\cM^\Z_{2m}$. We have used this relation in order to name elements of $tmf_*$ explicitly.

In the following we recall some material from \cite{rezk}.
Let $A$ be a commutative ring. A generalized Weierstrass equation over $A$
is an equation of the form
\begin{equation}\label{fhzwfewfewf87weff}
Y^2Z+a_1XYZ+a_3YZ^2=X^3+a_2X^2Z+a_4XZ^2+a_6Z^3
\end{equation} with $a_i\in A$.  We let 
$$W_0:=\spec  \Z[a_1,a_2,a_3,a_4,a_6]$$ be the affine scheme over $\Z$   representing 
Weierstrass equations.  
We furthermore let 
$$\xymatrix{\cC\ar@{^{(}->}[r]\ar[dr]& W_0\times \P^2\ar[d]^{pr_{W_0}}\\ &W_0} $$
 denote the universal Weierstrass curve defined by equation (\ref{fhzwfewfewf87weff}).

By $G\subseteq PGL_{3,\Z}$ we denote the subgroup scheme  
 of  automorphisms of $\P^2_\Z$ which respect
generalized Weierstrass equations. 
The group scheme $G$ is affine and  acts on the universal Weierstrass curve $\cC\to W_0$.

The scheme $W_1:=G\times W_0$ is the scheme of morphisms of the action groupoid $W_1\rightrightarrows W_0$  
 in the category of schemes. 
It will be called the Weierstrass groupoid and is the spectrum of the so-called Weierstrass Hopf algebroid
$(\cO_{W_0}(W_0),\cO_{W_1}(W_1))$. 

 The stack represented by the Weierstrass  groupoid will be denoted by $\WWW$.  The action of $G$ on the universal Weierstrass curve $\cC$ gives an  action groupoid
$G\times\cC\rightrightarrows  \cC$, and by $\CCC$ we denote the corresponding stack.
The $G$-equivariant projection $\cC\to W_0$ induces a morphism of stacks
$\pi:\CCC\to \WWW$.

For every ring $A$ the point $[0:1:0]\in \P^2_A$  belongs to each Weierstrass curve and is fixed by every automorphism in $G(A)$. We therefore get a section 
$$\xymatrix{\CCC\ar[d]^\pi\\ \WWW\ar@/^1cm/@{.>}[u]^s}\ .$$ The section $s$ lies in the smooth locus of $\pi$.
 The pull-back of the relative cotangent bundle
$s^*T^*(\CCC/\WWW)\to \WWW$ is thus a line bundle over $\WWW$.  
Its sheaf of sections will be denoted by $ \omega$.

We define an evenly graded Hopf algebroid $(A_*,\Gamma_*)$ whose underlying graded rings are given by 
$$A_*:=\bigoplus_{n\ge 0} \omega^{2n}(W_0)\ ,\quad \Gamma_*:= \bigoplus_{n\ge 0} \omega^{2n} (W_1)\ ,$$ 
using the map $W_1\stackrel{s}{\to} W_0\to \WWW$ in order to view
$W_1$ as a scheme over $\WWW$. The structure maps of this Hopf algebroid are induced in the canonical way.
The spectrum $tmf$ is characterized by the property that the complex $({}^YE_1^{*,2*}(tmf),d_1)$ has a preferred bi-graded isomorphism to the cobar complex $(C^*(A_*),d^{cobar})$ of the Hopf-algebroid $(A_*,\Gamma_*)$. This isomorphism will be used 
to name elements of ${}^YE_r^{*,*}(tmf_{(p)})$ explicitly.

Let $\bar \MMM$ denote the Deligne-Rapoport compactification  of the moduli stack of elliptic curves over $\Z$, i.e. the stack of generalized elliptic curves with irreducible geometric fibers, denoted $M_1$ in \cite{MR0337993}. Since generalized elliptic curves admit
Weierstrass equations (Zariski locally), there
is a unique (1-)morphism $i:\bar\MMM\to\WWW$
such that $i^*(\CCC)$ is isomorphic to the universal curve over $\bar\MMM$. One can check that $i$ is an open immersion.

The algebro-geometric definition of the space of 
integral modular forms of weight $m$ is
$${}^{alg}\cM_{m}^\Z:=H^0(\bar \MMM,i^* \omega^{m})\ .$$

Since $\WWW$ is normal and
the codimension of $\WWW\setminus i(\bar \MMM)$   is $\ge 2$ the inclusion $i:\bar \MMM \to \WWW$ induces an isomorphism
$$i^*:H^0( \WWW,\omega^{m})\to H^0(\bar\MMM,i^*\omega^{m})\ .$$
The following isomorphism will be used to identify elements in ${}^YE_2^{0,*}(tmf)$
with modular forms.
\begin{equation}\label{udiqdqwdqwdqwdqwd}
 {}^{alg}\cM_{2m}^\Z \stackrel{(i^*)^{-1}}{\to}
H^0(\WWW,\omega^{2m}) \cong 
H^0(C^*(A_* ,d^{cobar}) )_{2m}
\cong  {}^YE_2^{0,4m}(tmf)\ ,
\end{equation}
where $$H^0(C^*(A_* ,d^{cobar}) )=\bigoplus_{m\ge 0} H^0(C^*(A_* ,d^{cobar}) )_{2m}$$
is the decomposition of the cohomology by degree.

We now relate the algebraists'  version ${}^{alg}\cM^\Z_{2m}$ of modular forms with
our formal power series version  $\cM^\Z_{2m}\subset \Z[[q]]$.
 In homogeneous coordinates $x:=X/Z$, $y:=Y/Z$ we can consider the  $1$-form
$$\eta\in \Gamma(\cC, T^*(\cC/W_0))\ ,\quad  \eta:=\frac{dx}{2y+a_1x+a_3}=\frac{dy}{3x^2+2a_2x+a_4-a_1y}\ .$$
For all $n\ge 0$ the section $s^*\eta^n\in \omega^n(W_0)$ is nowhere vanishing and   
 induces isomorphisms of groups 
$$\eta^n:\cO_{W_0}(W_0)\stackrel{\sim}{\to} \omega^n(W_0)\ .$$

The homomorphism $Tate:\spec\Z[[q]]\to W_0$ defined by 
$$a_1\mapsto 1\ ,a_2\mapsto 0\ ,a_3\mapsto 0\ , a_4\mapsto B\ , a_6\mapsto  C$$
with $B,C\in \Z[[q]]$  given by

\[ B:=-5\sum\limits_{n\ge 1}n^3\cdot\frac{q^n}{1-q^n}, C:=-\frac{1}{12}\sum\limits_{n\ge 1}(7n^5+5n^3)\cdot\frac{q^n}{1-q^n}\]

defines the Tate curve over $\Z[[q]]$. 
It induces an injective map
$$r:{}^{alg}\cM_{2m}^\Z \stackrel{(i^*)^{-1}}{\to } H^0(\WWW,\omega^{2m}) \hookrightarrow
\omega^{2m}(W_0)\stackrel{(\eta^{2m})^{-2}}{\to}\cO_{W_0}(W_0)\stackrel{Tate^*}{\to} \Z[[q]]\ .$$
The following Lemma is a consequence of the $q$-expansion principle, \cite[Theorem 1.6.1]{katzpadic}.
\begin{lem}
The map $r$ defined above induces
an isomorphism $r:{}^{alg}\cM_{2m}^\Z\stackrel{\sim}{\to} \cM_{2m}^\Z$.
\end{lem}

If we combine the $q$-expansion principle with (\ref{udiqdqwdqwdqwdqwd})
we get the identification
\begin{equation}\label{ghagdastdtztzqw}
 {}^YE_2^{0,4m}(tmf)\cong \cM_{2m}^\Z
\end{equation}
used in the present paper.

We now turn to the structure of the groups $T_{2m}$.
For $\nu\ge 0$ we denote by $p_\nu :\R[[q]]\to\R$ the projection onto the
$\nu$-th coefficient.
Let $m\ge 2$ and  $k_m:=\dim_\C \cM^\C_{2m}$ be the dimension of the space of modular forms of weight $2m$.  
\begin{lem}\label{invariant}
\begin{itemize}
\item[i)] There exists a $\Z$-basis $f_0,\ldots, f_{k_m-1}\in \cM^\Z_{2m}$ such that $p_i(f_j)=\delta_{i,j}$ for $0\leq i,j\leq k_m-1$.
\item[ii)] For a basis as in part i), the map
\[ \alpha: T_{2m}=\frac{\R[[q]]}{\Z[[q]]+\cM^\R_{2m}}\longrightarrow \prod_{\nu\ge k_m}\left( \R/\Z\right)\ , \quad  \, f\mapsto \left(\left[ p_\nu(f-\sum\limits_{i=0}^{k-1} p_i(f)f_i)\right]\right)_{\nu\ge k_m}\]
is well-defined and an isomorphism of abelian groups.
\end{itemize}
\end{lem}
\proof It is easy to see that $i)$ implies $ii)$. 
We now show  $i)$.
There is an isomorphism $[ \bar \MMM]\cong \P^1_\Z$ of the coarse moduli space $[ \bar \MMM]$ of the  Deligne-Rapoport compactification of the moduli stack of elliptic
curves $\bar \MMM$, \cite[VI, Th\'eor\`eme 1.1]{MR0337993}. Let $\pi: \bar \MMM \to [\bar \MMM]$ denote the canonical projection. 
One can check that the sheaf $\pi_*i^*\omega^{2m}$ is locally free of rank one.
Hence there exists a unique function
$n:\nat_0\to \Z$ such that
$$pi_*i^*\omega^{2m}\simeq\cO(n(m)) \ .$$
The projection
$\pi:\bar \MMM\to  [\MMM]$ induces an isomorphism 
$$\pi^*:H^0(\P^1_\Z,\cO(n(m)))\stackrel{\sim}{\to}H^0(\bar \MMM,i^* \omega^{2m}) \ .$$
The map
$$H^0(\P^1_\Z,\cO(n(m)))\stackrel{\pi^*}{\to }H^0(\bar \MMM,i^* \omega^{2m})={}^{alg}\cM_{2m}^\Z\stackrel{r}{\to} \Z[[q]]$$ corresponds
to the Taylor expansion at a point $\infty\in \P^1_\Z$,
\cite[Chapter 8.11]{KatzMazur}.

 Its is now elementary to write down a basis of 
$H^0(\P^1_\Z,\cO(n(m)))$ whose Taylor expansion at $\infty$
has the property analogous to $i)$. 
\hB
 
The following corollary will be used later. Let $\Delta\in \cM^{\Z}_{12}$ be the unique cusp form normalized such that $\Delta=q+\dots$.
\begin{kor}\label{zdwdqwdqwdwddqwdq}
For every $k\ge 0$ and integer $a\in \nat$ the element
$\left[\frac{1}{a}\Delta^k \right]\in T_{12k+2}$
has order $a$.
\end{kor}
\proof 
The case $k=0$ is clear since $\cM_{2}=0$. Let us now assume that $k\ge 1$.
We have $\dim  \cM_{12k+2}=k$.  Since $p_\nu(\frac{1}{a}\Delta^k)=0$ for $\nu=0,\dots,k-1$
we have  
$$ \alpha(\left[ \frac{1}{a}\Delta^k\right])=\left(p_\nu(\frac{1}{a}\Delta^k)\right)_{\nu\ge k} =\left( \left[ \frac{1}{a} \right], \left[ \frac{a_1}{a}\right],\ldots\right)\in(\R/\Z)^\nat$$
for suitable $a_1,\ldots\in\Z$. 
 Therefore 
$\alpha(\left[ \frac{1}{a}\Delta^k\right])$ has order $a$. \hB

\section{Table of $MSpin_*$}\label{udiqdqwdwqdqwdwqdwqd}

Anderson-Brown-Peterson \cite{MR0190939} have obtained   the following additive decomposition $$MSpin^{\hspace{0.2em} \hat{\mbox{}}}_2\cong \prod_{n(J)\:\mbox{even}} ko^{\hspace{0.2em} \hat{\mbox{}}}_2\langle n(J)\rangle \vee \prod_{1\neq n(J)\:\mbox{odd}} ko^{\hspace{0.2em} \hat{\mbox{}}}_2\langle n(J)-2\rangle \vee  \prod_{i=0}^\infty \prod_{j=1}^{\dim(Z_i)} K(\Z/2\Z,i)\ $$
of the $2$-completion of $MSpin$.
In the first two  products the index $J$ runs over all unordered tuples
$J=(j_1,\dots,j_k)$ for $k\ge 0$ and $j_i\ge 2$, $n(J):=\sum_{i=1}^k j_i$, and $Z:=\bigoplus_{i=0}^\infty Z_i$ is some
$\nat$-graded $\F_2$-vector space.  Furthermore, $ko^{\hspace{0.2em} \hat{\mbox{}}}_2\langle k\rangle$ denotes the $2$-completed and $k$-connective cover of $ko$.
The Poincar\'e series given in  \cite[Thm.1.11]{MR0190939} can be combined to
determine the dimensions $\dim(Z_i)$ for all $i\ge 0$. Using MAPLE and these formulas we have compiled
the following table.
Recalling that there is no odd torsion in $MSpin_*$, this table determines 
$MSpin_*$ additively for $*\leq 127$. The final column separates off the
torsion which is accounted for by the $K(\Z/2\Z,i)$'s, i.e. those Spin-manifolds which are detected by Stiefel-Whitney classes. 
For the purpose of the present paper, note that the familiar structure
of $ko_*$ implies that $ MSpin_{4m-1}$ is entirely accounted for by the
$K(\Z/2\Z,i)$'s. In the smallest dimension $MSpin_{39}\simeq\Z/2\Z$
one can use computations of Mahowald/Gorbunov to show that $  MString_{39}\to
\ MSpin_{39}$ is zero, i.e. $A_{39}=MSpin_{39}$ but for all
$m\ge 11$ such that $ MSpin_{4m-1}\neq 0$ (the table strongly suggests
this is in fact true for {\em all} $m\ge 11$) we cannot decide whether 
or not $A_{4m-1}\subseteq MString_{4m-1}$ is an equality.

\newpage 
\thispagestyle{empty}
\begin{tabular}{|c|c|c|c|}\hline $i$&$\dim_\Q(MSpin\Q_i)$&$\dim_{\mathbb{F}_2}(MSpin_{i,tors})$& $ \dim_{\mathbb{F}_2} Z_i$ 
\\\hline 0&1&0&0
\\\hline 1&0&1&0
\\\hline 2&0&1&0
\\\hline 3&0&0&0
\\\hline 4&1&0&0
\\\hline 5&0&0&0
\\\hline 6&0&0&0
\\\hline 7&0&0&0
\\\hline 8&2&0&0
\\\hline 9&0&2&0
\\\hline10&0&3&0
\\\hline11&0&0&0
\\\hline12&3&0&0
\\\hline13&0&0&0
\\\hline14&0&0&0
\\\hline15&0&0&0
\\\hline16&5&0&0
\\\hline17&0&5&0
\\\hline18&0&7&0
\\\hline19&0&0&0
\\\hline20&7&1&1
\\\hline21&0&0&0
\\\hline22&0&1&1
\\\hline23&0&0&0
\\\hline24&11&0&0
\\\hline25&0&11&0
\\\hline26&0&15&0
\\\hline27&0&0&0
\\\hline28&15&2&2 
\\\hline29&0&1&1
\\\hline30&0&3&3
\\\hline31&0&0&0
\\\hline32&22&1&1
\\\hline33&0&23&1
\\\hline34&0&31&1
\\\hline35&0&0&0
\\\hline36&30&6&6
\\\hline37&0&2&2
\\\hline38&0&7&7
\\\hline39&0&1&1
\\\hline40&42&4&4
\\\hline41&0&45&3
\\\hline42&0&60&4
\\\hline43&0&2&2
\\\hline
\end{tabular}

\newpage 
\thispagestyle{empty}
\begin{tabular}{|c|c|c|c|}\hline $i$&$\dim_\Q(MSpin\Q_i)$&$\dim_{\mathbb{F}_2}(MSpin_{i,tors})$& $ \dim_{\mathbb{F}_2} Z_i$ 
  \\\hline   44& 56& 14& 14
                    \\\hline            45& 0& 6& 6
                    \\\hline            46& 0& 17& 17
                    \\\hline             47& 0& 4& 4
                    \\\hline           48& 77& 11& 11
                    \\\hline            49& 0& 86& 9
                    \\\hline           50& 0& 114& 13
                    \\\hline             51& 0& 7& 7
                    \\\hline           52& 101& 31& 31
                    \\\hline            53& 0& 15& 15
                    \\\hline            54& 0& 38& 38
                    \\\hline            55& 0& 13& 13
                    \\\hline           56& 135& 29& 29
                    \\\hline           57& 0& 159& 24
                    \\\hline           58& 0& 210& 34
                    \\\hline            59& 0& 22& 22
                    \\\hline           60& 176& 67& 67
                    \\\hline            61& 0& 38& 38
                    \\\hline            62& 0& 80& 80
                    \\\hline            63& 0& 36& 36
                    \\\hline           64& 231& 70& 70
                    \\\hline           65& 0& 290& 59
                    \\\hline           66& 0& 379& 82
                    \\\hline            67& 0& 58& 58
                    \\\hline          68& 297& 142& 142
                    \\\hline            69& 0& 90& 90
                    \\\hline           70& 0& 169& 169
                    \\\hline            71& 0& 92& 92
                    \\\hline          72& 385& 158& 158
                    \\\hline           73& 0& 521& 136
                    \\\hline           74& 0& 676& 186
                    \\\hline           75& 0& 143& 143
                    \\\hline          76& 490& 291& 291
                    \\\hline           77& 0& 205& 205
                    \\\hline           78& 0& 347& 347
                    \\\hline           79& 0& 219& 219
                    \\\hline          80& 627& 343& 343
                    \\\hline           81& 0& 931& 304
                    \\\hline          82& 0& 1196& 404
                    \\\hline           83& 0& 330& 330
                    \\\hline          84& 792& 589& 589
                    \\\hline           85& 0& 448& 448
                    \\\hline           86& 0& 698& 698
                    \\\hline           87& 0& 494& 494
\\\hline
\end{tabular}

\newpage 
\thispagestyle{empty}
\begin{tabular}{|c|c|c|c|}\hline $i$&$\dim_\Q(MSpin\Q_i)$&$\dim_{\mathbb{F}_2}(MSpin_{i,tors})$& $ \dim_{\mathbb{F}_2} Z_i$
  \\\hline         88& 1002& 721& 721
                    \\\hline          89& 0& 1658& 656
                    \\\hline          90& 0& 2103& 848
                    \\\hline           91& 0& 729& 729
                    \\\hline        92& 1255& 1171& 1171
                    \\\hline           93& 0& 952& 952
                    \\\hline          94& 0& 1385& 1385
                    \\\hline          95& 0& 1068& 1068
                    \\\hline        96& 1575& 1472& 1472
                    \\\hline          97& 0& 2948& 1373
                    \\\hline          98& 0& 3689& 1731
                    \\\hline          99& 0& 1550& 1550
                    \\\hline        100& 1958& 2296& 2296
                    \\\hline         101& 0& 1967& 1967
                    \\\hline         102& 0& 2706& 2706
                    \\\hline         103& 0& 2233& 2233
                    \\\hline        104& 2436& 2941& 2941
                    \\\hline         105& 0& 5239& 2803
                    \\\hline         106& 0& 6461& 3451
                    \\\hline         107& 0& 3194& 3194
                    \\\hline        108& 3010& 4438& 4438
                    \\\hline         109& 0& 3969& 3969
                    \\\hline         110& 0& 5215& 5215
                    \\\hline         111& 0& 4539& 4539
                    \\\hline        112& 3718& 5760& 5760
                    \\\hline        113& 0& 9312& 5594
                    \\\hline         114& 0& 11311& 6746
                    \\\hline         115& 0& 6411& 6411
                    \\\hline        116& 4565& 8470& 8470
                    \\\hline         117& 0& 7839& 7839
                    \\\hline         118& 0& 9925& 9925
                    \\\hline         119& 0& 9005& 9005
                    \\\hline       120& 5604& 11086& 11086
                    \\\hline        121& 0& 16544& 10940
                    \\\hline        122& 0& 19796& 12954
                    \\\hline        123& 0& 12582& 12582
                    \\\hline       124& 6842& 15963& 15963
                    \\\hline        125& 0& 15193& 15193
                    \\\hline        126& 0& 18656& 18656
                    \\\hline        127& 0& 17493& 17493

\\\hline
\end{tabular}

\end{document}